      \documentclass[11pt]{amsart}
      \usepackage[mathscr]{eucal}
      \usepackage{amsmath,amsfonts}
\def\mapright#1{\smash{
\mathop{\rg}\limits^{#1}}}
\def\mapdown#1{\bigg\downarrow
\rlap{$\vcenter{\hbox{$\scriptstyle#1$}}$}}

\def\rg{\hbox to 30pt{\rightarrowfill}}
\def\lg{\hbox to 30pt{\leftarrowfill}}

      \parskip\smallskipamount
          \newtheorem{theorem}{Theorem}[section]
      
      \newtheorem{proposition}[theorem]{Proposition}
      \newtheorem{corollary}[theorem]{Corollary}
      \newtheorem{lemma}[theorem]{Lemma}
      \newtheorem{example}[theorem]{Example}

      \newtheorem{remark}[theorem]{Remark}
      \makeatletter
      \@addtoreset{equation}{section}
      \makeatother

\newcommand{\HH}{{\mathbb H}}
      \newcommand{\BB}{{\mathbb B}}
      \newcommand{\CC}{{\mathbb C}}
      \newcommand{\NN}{{\mathbb N}}

      \newcommand{\DD}{{\mathbb D}}
      
      \newcommand{\FF}{{\mathbb F}}

      \newcommand{\cA}{{\mathcal A}}
      
      \newcommand{\cC}{{\mathcal C}}
      \newcommand{\cD}{{\mathcal D}}

      \newcommand{\cE}{{\mathcal E}}
      
      \newcommand{\cG}{{\mathcal G}}
      \newcommand{\cH}{{\mathcal H}}
      \newcommand{\cK}{{\mathcal K}}
      \newcommand{\cL}{{\mathcal L}}
      \newcommand{\cM}{{\mathcal M}}
      \newcommand{\cN}{{\mathcal N}}
      \newcommand{\cP}{{\mathcal P}}
      
      \newcommand{\cS}{{\mathcal S}}

      \newcommand{\cY}{{\mathcal Y}}
      
      \newcommand{\cW}{{\mathcal W}}
      \newcommand{\cZ}{{\mathcal Z}}

      \newcommand{\rank}{\hbox{\rm{rank}}\,}

      \newdimen\expt
      \def\boxit#1{\setbox0\hbox{$\displaystyle{#1}$}
            \hbox{\lower.4\expt
       \hbox{\lower3\expt\hbox{\lower\dp0
            \hbox{\vbox{\hrule height.4\expt
       \hbox{\vrule width.4\expt\hskip3\expt
            \vbox{\vskip3\expt\box0\vskip2\expt}%
       \hskip3\expt\vrule width.4\expt}\hrule height.4\expt}}}}}}
      
\begin{document}
       \pagestyle{myheadings}
      \markboth{ Gelu Popescu}{Operator theory on noncommutative varieties}

      \title [  Operator theory on noncommutative varieties ]
      { Operator theory on noncommutative varieties}
        \author{Gelu Popescu}
      \thanks{Research supported in part by an NSF grant}

      \address{Department of Mathematics, The University of Texas
      at San Antonio \\ San Antonio, TX 78249, USA}
      \email{\tt gelu.popescu@utsa.edu}

      \begin{abstract}

We develop a dilation theory for row contractions $T:=[T_1,\ldots, T_n]$
subject to constraints such as
$$
p(T_1,\ldots, T_n)=0,\quad p\in \cP,
$$
where $\cP$ is a set of noncommutative polynomials. The model $n$-tuple is the universal row contraction $[B_1,\ldots, B_n]$ satisfying the same constraints as $T$, which turns out to be, in a certain sense, the  {\it maximal  constrained piece} of the $n$-tuple $[S_1,\ldots, S_n]$ of left creation operators on the full Fock space on $n$ generators. The theory is based on a class of
{\it constrained Poisson kernels} associated with $T$  and representations of the $C^*$-algebra generated  by $B_1,\ldots, B_n$ and the identity.
Under natural conditions  on the constraints we have uniqueness for the minimal dilation.

A characteristic function $\Theta_T$ is associated with any (constrained) row contraction $T$ and it is proved that
$$I-\Theta_T \Theta_T^*=K_TK_T^*,
$$
where $K_T$ is the  (constrained) Poisson kernel of $T$.
Consequently, for {\it pure constrained} row contractions, we show that
the characteristic function  is  a  complete unitary invariant and provide a   model. We show that the curvature invariant and Euler characteristic
asssociated with a Hilbert module generated by  an arbitrary
(resp.~commuting) row contraction $T$ can be expressed only in terms of the
(resp.~constrained) characteristic function of $T$.

We provide a commutant lifting theorem for pure constrained row contractions and obtain a Nevanlinna-Pick  interpolation result  in our setting.

 \end{abstract}
      \maketitle
      \section*{Introduction}

    In  the last two decades, there has been exciting progress in  multivariable   dilation theory, in the attempt to extend the classical Nagy-Foia\c s theory of contractions \cite{SzF-book}.  In the noncommutative case, significant results were obtained in \cite{Fr}, \cite{Bu}, \cite{Po-models}, \cite{Po-isometric},
\cite{Po-charact},
 \cite{Po-intert}, and recently in \cite{DKS}. Some of these results were further extended by Muhly and Solel \cite{MuSo} to  representations of tensor algebras over $C^*$-correspondences.

The present paper develops a dilation theory for row contractions  subject to constraints   determined by   sets of noncommutative polynomials.  Our theory includes, in particular, the commutative (see \cite{Dr}, \cite{Po-poisson}, and  \cite{Arv}) and $q$-commutative (see \cite{APo}, \cite{BB})
cases, while the standard noncommutative dilation theory
for row contractions serves as a ``universal model".

An $n$-tuple $T:=[T_1,\ldots, T_n]$ of  bounded linear operators acting on a common Hilbert space $\cH$ is called  {\it row contraction} if
$$
TT^*=T_1T_1^*+\cdots+T_nT_n^*\leq I.
$$
A distinguished role among row contractions is played by
the $n$-tuple $S:=[S_1,\ldots, S_n]$ of {\it left creation operators} on the full Fock space $F^2(H_n)$, which satisfies the noncommutative von Neumann inequality \cite{Po-von}
$$\|p(T_1,\ldots, T_n)\|\leq \|p(S_1,\ldots, S_n)\|
$$
for any polynomial $p(X_1,\ldots, X_n)$ in $n$ noncommuting indeterminates. For the classical von Neumann inequality \cite{vN} (case $n=1$) and a nice survey,  we refer to Pisier's book \cite{Pi}.
Based on the left creation operators and their representations, a noncommutative dilation theory and model theory for row contractions was developed in \cite{Fr}, \cite{Bu}, \cite{Po-models}, \cite{Po-isometric},
\cite{Po-charact},
 \cite{Po-intert}, etc.
Assume now that $T$ is subject to the constraints
$$
T_iT_j=T_jT_i,\quad i,j=1,\ldots,n.
$$
In this commutative case, the noncommutative dilation theory can be applied but, in many respects, it is not satisfactory due to the fact that the model shift $S:=[S_1,\ldots, S_n ]$ does not satisfy the same constraints as $T$. However, the universal commutative row contraction  is a
{\it piece} of  $S$, namely
$$
B:=[B_1,\ldots, B_n],\qquad B_i:=P_{F_s^2} S_i|F^2_s,\quad i=1,\ldots, n,
$$
where $F_s^2\subset F^2(H_n)$ is the symmetric  Fock space.
In this setting, the natural von Neumann inequality (see \cite{Dr}, \cite{Po-poisson}, and  \cite{Arv})
is
$$\|p(T_1,\ldots, T_n)\|\leq \|p(B_1,\ldots, B_n)\|
$$
for any polynomial \, $p(z_1,\ldots, z_n)$\,  in $n$ commuting indeterminates.
A dilation theory for commuting row contractions based on the model shift $B:=[B_1,\ldots, B_n]$ and its representations was considered by Drury \cite{Dr} and the author \cite{Po-poisson} to a certain extent, and by Arveson \cite{Arv} in greater details.
This circle of ideas was extended  to row contractions satisfying
the constraints
$$T_jT_i=q_{ij} T_iT_j,\quad 1\leq i< j\leq n,
$$
where $q_{ij}\in \CC$. In this setting, a von Neumann inequality was obtained  by Arias and the author in \cite{APo}.
 This was used further by B.V.R.~Bhat and T.~Bhattacharyya \cite{BB}
to obtain   a model theory  for $q$-commuting row contractions.

In this paper, we develop a dilation theory for row contractions $T:=[T_1,\ldots, T_n]$
subject to  more general constraints such as
$$
p(T_1,\ldots, T_n)=0,\quad p\in \cP,
$$
where $\cP$ is a set of noncommutative polynomials.
If $T$ is a {\it pure} row contraction, then $\cP$ can be any WOT-closed two-sided ideal of the noncommutative analytic Toeplitz algebra $F_n^\infty$.
The model $n$-tuple is the universal row contraction $[B_1,\ldots, B_n]$ satisfying the same constraints as $T$, which turns out to be, in a certain sense, the  {\it  maximal  constrained piece} of the $n$-tuple $[S_1,\ldots, S_n]$ of left creation operators on the full Fock space with $n$ generators.

In Section \ref{Beurling}, we provide basic results concerning the {\it constrained  shift }  $[B_1,\ldots, B_n]$ and  the $w^*$-closed algebra (resp.~$C^*$-algebra) generated by $B_1,\ldots, B_n$ and the identity. We  obtain a Beurling type theorem characterizing the invariant subspaces under each operator \linebreak $B_1\otimes I_\cH,\ldots, B_n\otimes I_\cH$, and a characterization of cyclic co-invariant subspaces under the same operators.
We also provide Wold type decompositions for  nondegenerate $*$-representations of the $C^*$-algebra $C^*(B_1,\ldots, B_n)$ and prove that two constrained shifts $[B_1\otimes I_\cH,\ldots, B_n\otimes I_\cH]$ and
$[B_1\otimes I_\cK,\ldots, B_n\otimes I_\cK]$ are similar if and only if
$\dim \cH=\dim \cK$.

In Section \ref{Dilations}, we develop a dilation theory
for constrained row contractions.
The theory is based on a class of
{\it constrained Poisson kernels} (see \cite{Po-poisson}, \cite{APo},  \cite{Po-unitary}, and \cite{BC} for $n=1$) associated with
$T:=[T_1,\ldots, T_n]$  and representations of the $C^*$-algebra generated  by $B_1,\ldots, B_n$ and the identity.
In particular, if the  set $\cP$ consists of homogenous polynomials, then we show that
 there exists a  Hilbert space $\cK_\pi$ such that $\cH$ can be identified with a subspace of \linebreak
$\tilde\cK:=(\cN_J\otimes \overline{\Delta_T\cH})\oplus \cK_\pi$
and
$$
T_i^*=V_i^*|\cH,\quad i=1,\ldots, n,
$$
where  $\Delta_T:=(I-T_1T_1^*-\cdots -T_nT_n^*)^{1/2}$,
$$
V_i:=\left[\begin{matrix}
B_i\otimes I_{\overline{\Delta_T\cH}}&0\\0&\pi(B_i)
\end{matrix}\right],\quad i=1,\ldots,n,
$$
and  $\pi:C^*(B_1,\ldots, B_n)\to B(\cK_\pi)$ is a $*$-representation
which annihilates the compact operators and
$$
\pi(B_1)\pi(B_1)^*+\cdots +\pi(B_n)\pi(B_n)^*=I_{ \cK_\pi}.
$$
Under certain  natural conditions  on the constraints, we have uniqueness for the minimal dilation of $T$. We introduce  and evaluate
the {\it dilation index}, a numerical invariant for row contractions,  and show that it does not depend on the constraints.

In Section \ref{Characteristic} and Section \ref{Characteristic-constr},  we provide new properties for the {\it standard characteristic function} $\Theta_T$ associated with an arbitrary row contraction $T$ (see \cite{Po-charact}), and  introduce  a new characteristic function  associated  with  constrained row contractions.
The {\it constrained characteristic function} is
   a   multi-analytic operator (with respect to the constrained shifts $B_1,\ldots, B_n$)
$$
\Theta_{J,T}:\cN_J\otimes \cD_{T^*}\to \cN_J\otimes \cD_T
$$
uniquely defined by the formal Fourier representation
$$ -I_{\cN_J}\otimes T+
\left(I_{\cN_J}\otimes \Delta_T\right)\left(I_{{\cN_J}\otimes \cH}-\sum_{i=1}^n W_i\otimes T_i^*\right)^{-1}\\
\left[W_1\otimes I_\cH,\ldots, W_n\otimes I_\cH
\right] \left(I_{\cN_J}\otimes \Delta_{T^*}\right)
$$
(see Section \ref{Characteristic-constr} for notations).
We prove a factorization result  for
the   constrained  characteristic function, namely,
$$I-\Theta_{J,T} \Theta_{J,T}^*=K_{J,T} K_{J,T}^*,
$$
where $K_{J,T}$ is the  constrained  Poisson kernel
associated with $T$.
Consequently, for the class of pure constrained row contractions,
we show that
the characteristic function  is  a complete unitary invariant and provide a model. All the results of Section
\ref{Characteristic-constr} apply, in particular, to commutative row contractions.

In Section \ref{Commutant}, we obtain a commutant lifting theorem for pure constrained row contractions and a Nevanlinna-Pick  interpolation result  in our setting. These results are based on the more general noncommutative commutant lifting theorem  (see \cite{Po-isometric}, \cite{Po-intert}) and some  results from previous sections.

The above-mentioned factorization result
for
the   characteristic function
has important consequences in multivariable operator theory.  We point out some of them which are considered in  Section \ref{Characteristic} and Section \ref{Characteristic-constr}.
 In \cite{Arv2}, Arveson
    introduced  a notion of
        curvature   and Euler characteristic  for finite rank
       contractive Hilbert modules over $\CC[z_1,\ldots, z_n]$,
       the complex unital algebra of all polynomials in $n$ commuting variables.
  The canonical operators  $T_1,\dots, T_n $ associated with  the $\CC[z_1,\dots, z_n]$-module structure are commuting and $T:=[T_1,\ldots, T_n]$ is a row contraction
with $\rank \Delta_T<\infty$.
    Noncommutative analogues of these notions were introduced  and studied
    by the author
   in  \cite{Po-curvature} and, independently, by D.~Kribs  \cite{Kr}.
   In this paper, we show that the curvature and the  Euler characteristic (in both  the commutative and noncommutative  case) depend only on  the properties of the characteristic function of $T$.

For example, in the commutative case,
if  $T:=[T_1,\ldots, T_n]$  is  a commutative  row contraction  with $\rank \Delta_T<\infty$, and   $K(T)$ and $\chi(T)$ denote Arveson's curvature and Euler characteristic, respectively,
then we prove  that
\begin{equation*}\begin{split}
K(T)&=\int_{\partial \BB_n}\lim_{r\to 1}\text{\rm trace}\, [I-\Theta_{J_c,T}(r\xi)\Theta_{J_c,T}(r\xi)^*] d\sigma(\xi)\\
&=
\rank \Delta_T-(n-1)! \lim_{m\to \infty}
\frac {\text{\rm trace}\, [\Theta_{J_c,T}\Theta_{J_c,T}^* (Q_m\otimes I_{\cD_T})]}
{n^m},
\end{split}\end{equation*}
where $Q_m$ is the projection of  Arveson's space $H^2$ onto the subspace   of homogeneous polynomials of degree $m$,
and
$$
\chi(T)=n!
\lim_{m\to \infty}
\frac {\rank\left[\left(I-\Theta_{J_c,T}\Theta_{J_c,T}^*\right) \left(Q_{\leq m}\otimes I_{\cD_T}\right)\right]} {m^n},
$$
where $Q_{\leq m}$ is the projection of $H^2$ onto the subspace  of all polynomials of degree $\leq m$.
Here, the operator
$\Theta_{J_c,T}:H^2\otimes \cD_{T^*}\to H^2\otimes \cD_T$
stands for the constrained characteristic function associated with $T$, which, in this particular case,  is a
 multiplier  (multiplication operator)
defined by its symbol (for which we use the same notation)
 $$
\Theta_{J_c,T}(z):=
-T+\Delta_T(I-z_1T_1^*-\cdots -z_nT_n^*)^{-1} [z_1I_\cH,\ldots, z_nI_\cH]\Delta_{T^*},\quad z\in \BB_n,
$$
which is    a bounded  operator-valued analytic function on the
open unit ball  of $\CC^n$.

\smallskip

\bigskip

\section{Constrained shifts: invariant subspaces and Wold decompositions}\label{Beurling}

In this section we provide basic results concerning the
 {\it constrained  shift }  $[B_1,\ldots, B_n]$
associated with every WOT-closed two-sided ideal $J$ of the noncommutative analytic Toeplitz algebra $F_n^\infty$.
We  obtain a Beurling type theorem characterizing the invariant subspaces under each operator $B_1\otimes I_\cH,\ldots, B_n\otimes I_\cH$, and a characterization of cyclic co-invariant subspaces under the same operators.
 We also provide Wold type decompositions for the nondegenerate $*$-representations of the  Toeplitz $C^*$-algebra $C^*(B_1,\ldots, B_n)$   generated by $B_1,\ldots, B_n$ and the identity.
The dilation theory  developed in the next sections will be  based  on the
 constrained shift $[B_1,\ldots, B_n]$.

      Let $H_n$ be an $n$-dimensional complex  Hilbert space with orthonormal
      basis
      $e_1$, $e_2$, $\dots,e_n$, where $n\in\{1,2,\dots\}$.        We consider the full Fock space  of $H_n$ defined by
      $$F^2(H_n):=\bigoplus_{k\geq 0} H_n^{\otimes k},$$
      where $H_n^{\otimes 0}:=\CC 1$ and $H_n^{\otimes k}$ is the (Hilbert)
      tensor product of $k$ copies of $H_n$.
      Define the left creation
      operators $S_i:F^2(H_n)\to F^2(H_n), \  i=1,\dots, n$,  by
      $$
       S_i\varphi:=e_i\otimes\varphi, \quad  \varphi\in F^2(H_n).
      $$

      The    noncommutative analytic Toeplitz algebra   $F_n^\infty$
        and  its norm closed version,
        the noncommutative disc
       algebra  $\cA_n$,  were introduced by the author   \cite{Po-von}, \cite{Po-funct}, \cite{Po-disc} in
      connection
         with a multivariable noncommutative von Neumann inequality.
      $F_n^\infty$  is the algebra of left multipliers of
      the Fock space $F^2(H_n)$  and  can be identified with
       the
        weakly closed  (or $w^*$-closed) algebra generated by the left creation
      operators
         $S_1,\dots, S_n$  acting on   $F^2(H_n)$,
          and the identity.
           The noncommutative disc algebra $\cA_n$ is
          the  norm closed algebra generated by
         $S_1,\dots, S_n$,
          and the identity.
           When $n=1$, $F_1^\infty$
         can be identified
         with $H^\infty(\DD)$, the algebra of bounded analytic functions
          on the open unit disc. The  noncommutative analytic Toeplitz algebra $F_n^\infty$ can be viewed as a
           multivariable noncommutative
          analogue of $H^\infty(\DD)$.
       There are many analogies with the invariant
        subspaces of the unilateral
       shift on $H^2(\DD)$, inner-outer factorizations,
        analytic operators, Toeplitz operators, $H^\infty(\DD)$--functional
         calculus, bounded (resp.~spectral) interpolation, etc.
      %
%

Let $\FF_n^+$ be the unital free semigroup on $n$ generators
      $g_1,\dots,g_n$, and the identity $g_0$.
      The length of $\alpha\in\FF_n^+$ is defined by
      $|\alpha|:=k$ if $\alpha=g_{i_1}g_{i_2}\cdots g_{i_k}$, and
      $|\alpha|:=0$ if $\alpha=g_0$.
       If $T_1,\dots,T_n\in B(\cH)$, the algebra of all bounded operators on a Hilbert space $\cH$,  define \linebreak
      $T_\alpha :=  T_{i_1}T_{i_2}\cdots T_{i_k}$
      if $\alpha=g_{i_1}g_{i_2}\cdots g_{i_k}$, and
      $T_{g_0}:=I_\cH$.  Similarly, we denote $e_\alpha:=
e_{i_1}\otimes\cdots \otimes  e_{i_k}$ and $e_{g_0}:=1$.

      We need to recall from  \cite{Po-charact},
      \cite{Po-multi}, \cite{Po-von},  \cite{Po-funct},  \cite{Po-analytic}, and \cite{Po-central}
       a few facts
       concerning multi-analytic   operators on Fock
      spaces.
         We say that
       a bounded linear
        operator
      $M$ acting from $F^2(H_n)\otimes \cK$ to $ F^2(H_n)\otimes \cK'$ is
       multi-analytic
      if
      \begin{equation*}
      M(S_i\otimes I_\cK)= (S_i\otimes I_{\cK'}) M\quad
      \text{\rm for any }\ i=1,\dots, n.
      \end{equation*}
      Notice that $M$ is uniquely determined by   the ``coefficients''
        $\theta_{(\alpha)}\in B(\cK, \cK')$  given by
       $$
      \left< \theta_{(\tilde\alpha)}k,k'\right>:=  \left< M(1\otimes k), e_\alpha
      \otimes k'\right>,\quad
      k\in \cK,\ k'\in \cK',\ \alpha\in \FF_n^+,
      $$
      where $\widetilde\alpha$ is the reverse of $\alpha$, i.e.,
      $\widetilde\alpha=
      g_{i_k}\cdots g_{i_1}$ if
      $\alpha= g_{i_1}\cdots g_{i_k}$.
       We can associate with $M$ a unique formal Fourier expansion
      \begin{equation*}       M(R_1,\ldots, R_n):= \sum_{\alpha \in \FF_n^+} R_\alpha \otimes \theta_{(\alpha)}, \end{equation*}
      where $R_i:=U^* S_i U$, \ $i=1,\ldots, n$, are the right creation
      operators
      on $F^2(H_n)$ and
      $U$ is the (flipping) unitary operator on $F^2(H_n)$ mapping
        $e_{i_1}\otimes e_{i_2}\otimes\cdots\otimes e_{i_k}$ into
       $e_{i_k}\otimes\cdots\otimes e_{i_2}\otimes e_{i_1} $.
       Since  the operator $M$ acts like its Fourier representation on ``polynomials'',
  we will identify them for simplicity.
      Based on the noncommutative von Neumann
         inequality,         we proved that
        $$M(R_1,\ldots, R_n)=\text{\rm SOT}-\lim_{r\to 1}\sum_{k=0}^\infty
      \sum_{|\alpha|=k}
         r^{|\alpha|} R_\alpha\otimes \theta_{(\alpha)},
         $$
         where, for each $r\in [0,1)$, the series converges in the uniform norm.
      Moreover, the set of  all multi-analytic operators in
      $B(F^2(H_n)\otimes \cK,
      F^2(H_n)\otimes \cK')$  coincides  with
      $R_n^\infty\bar \otimes B(\cK,\cK')$,
      the WOT-closed operator space generated by the spatial tensor product, where
      $R_n^\infty=U^* F_n^\infty U$.
A multi-analytic operator is called inner if it is an isometry.

Let $J$ be a WOT-closed two-sided ideal of $F_n^\infty$
such that $J\neq F_n^\infty$, and define the subspaces of  the full Fock space $F^2(H_n)$ by setting
$$
\cM_J:=\overline{JF^2(H_n)}\quad \text{and}\quad \cN_J:=F^2(H_n)\ominus \cM_J.
$$
Notice that
$$
\cM_J=\overline{\{\varphi(1):\ \varphi\in
J\}}
\quad \text{and}\quad \cN_J=\bigcap_{\varphi\in J} \ker \varphi^*.
$$

Based on a Beurling type theorem \cite{Po-charact}   for the left creation  operators $S_1,\ldots, S_n$, a characterization of all $WOT$-closed two-sided ideal of $F_n^\infty$ was obtained  by Davidson and Pitts in \cite{DP2}.

One can easily obtain the following result.

\begin{lemma}\label{inv.sub}
Let $J$ be a WOT-closed two-sided ideal of $F_n^\infty$.
\begin{enumerate}
\item[(i)]
 If
$f(0)=0$ for any $f\in J$, then $1\in \cN_J$.
\item[(ii)]$\cN_J\neq 0$ if and only if $ J\neq F_n^\infty$.
\item[(iii)] The subspaces $\cN_J$ and $U\cN_J$ are invariant under $S_1^*,\ldots, S_n^*$ and  $R_1^*,\ldots, R_n^*$.
\end{enumerate}
\end{lemma}
\begin{proof}
The first part is obvious. The   second part  is a consequence of the fact  (see \cite{APo}, \cite{DP}) that, for any
$\varphi\in F_n^\infty$,
$$
d(\varphi, J)=\|P_{\cN_J}\varphi
(S_1,\ldots, S_n)|\cN_J\|.
$$
Part (iii) is straightforward.
\end{proof}

Define the {\it constrained  left} (resp.~{\it right}) {\it creation operators} by setting
$$B_i:=P_{\cN_J} S_i|\cN_J \quad \text{and}\quad W_i:=P_{\cN_J} R_i|\cN_J,\quad i=1,\ldots, n.
$$
Let $\cW(B_1,\ldots, B_n)$ be the $w^*$-closed algebra generated by
$B_1,\ldots, B_n$ and the identity. We proved in \cite{APo}
that $\cW(B_1,\ldots, B_n)$ has the $\mathbb A_1(1)$ property and therefore the $w^*$ and WOT topologies coincide on this algebra. Moreover,
we showed that
$$
\cW(B_1,\ldots, B_n)=P_{\cN_J}F_n^\infty |\cN_J=\{f(B_1,\ldots, B_n):\ f(S_1,\ldots, S_n)\in F_n^\infty\},
$$
where,  according to the $F_n^\infty$-functional calculus for row contractions \cite{Po-funct},
$$f(B_1,\ldots, B_n)=\text{\rm SOT-}\lim\limits_{r\to 1}f(rB_1,\ldots, rB_n).
$$
Note that if $\varphi\in J$, then $\varphi(B_1,\ldots, B_n)=0$.
Similar results hold for $\cW(W_1,\ldots, W_n)$, the $w^*$-closed algebra generated by
$W_1,\ldots, W_n$ and the identity. Moreover, we  proved in \cite{APo} that
$$
\cW(B_1,\ldots, B_n)^\prime=\cW(W_1,\ldots, W_n)\ \text{ and } \ \cW(W_1,\ldots, W_n)^\prime=\cW(B_1,\ldots, B_n),
$$
where $^\prime$ stands for the commutant.
An operator $M\in B(\cN_J\otimes \cK,\cN_J\otimes \cK')$ is called multi-analytic with respect to
the constrained shifts $B_1,\ldots, B_n$ if
$$
M(B_i\otimes I_{\cK})=(B_i\otimes I_{\cK'})M,\quad i=1,\ldots, n.
$$
If in addition $M$ is partially isometric, then we call it inner.
We recall from \cite{Po-central} that the set of all multi-analytic  operators with respect to
 $B_1,\ldots, B_n$ coincides  with
$$
\cW(W_1,\ldots, W_n)\bar\otimes  B(\cK,\cK')=P_{\cN_J\otimes \cK'}[R_n^\infty\bar\otimes B(\cK,\cK')]|\cN_J\otimes \cK,
$$
and a similar result holds for  the algebra $\cW(B_1,\ldots, B_n)$.

The next result provides  a Beurling type theorem characterizing
the invariant subspaces under the constrained shifts $B_1,\ldots, B_n$.

\begin{theorem}\label{Beur} Let $J\neq F_n^\infty$ be a WOT-closed two-sided ideal of $F_n^\infty$ and let $B_1,\ldots, B_n$ be the corresponding constrained left creation operators on $\cN_J$.
A subspace $\cM\subseteq \cN_J\otimes \cK$ is invariant under each
operator $B_i\otimes I_\cK$, \ $i=1,\ldots, n$, if and only if there exists
a Hilbert space  $\cG$ and an inner operator
$$\Theta(W_1,\ldots, W_n)\in \cW(W_1, \ldots, W_n)\bar\otimes B(\cG,\cK)$$ such that
\begin{equation}\label{Beu}
\cM=\Theta(W_1,\ldots, W_n)\left( \cN_J\otimes \cG\right).
\end{equation}
\end{theorem}

\begin{proof}
According to Lemma \ref{inv.sub}, the subspace $\cN_J\otimes K$ is invariant
under each operator $S_i^*\otimes I_\cK$, $i=1,\ldots, n$, and
$$(S_i^*\otimes I_\cK)|\cN_J\otimes \cK=B_i^*\otimes I_\cK, \quad i=1,\ldots,n. $$
Since the subspace $[\cN_J\otimes \cK]\ominus \cM$  is invariant under
 $B_i^*\otimes I_\cK$, \ $i=1,\ldots, n$, we deduce that it is also invariant under  each  operator $S_i^*\otimes I_\cK$, \ $i=1,\ldots, n$.
Therefore, the subspace
\begin{equation}\label{E}
\cE:=[F^2(H_n)\otimes \cK]\ominus\{[\cN_J\otimes \cK]\ominus \cM\}=
[\cM_J\otimes \cK]\oplus \cM
\end{equation}
is invariant under $S_i\otimes I_\cK$, \ $i=1,\ldots, n$, where $M_J:=F^2(H_n)\ominus \cN_J$.
Using the Beurling type characterization
of the invariant subspaces under the left creation operators (see
Theorem 2.2 from \cite{Po-charact}) and the characterization of multi-analytic operators from \cite{Po-analytic} (see also \cite{Po-central}), we find a Hilbert space $\cG$ and  an inner multi-analytic operator
$$\Theta(R_1,\ldots,R_n)\in R_n^\infty\bar \otimes B(\cG,\cK)$$
such that
$$
\cE=\Theta(R_1,\ldots,R_n)[F^2(H_n)\otimes \cG],
$$
where $\Theta(R_1,\ldots, R_n)$ is essentially unique up to a unitary diagonal multi-analytic operator.
Since
$\Theta(R_1,\ldots,R_n)$ is an isometry, we have
\begin{equation}\label{proj}
P_\cE=\Theta(R_1,\ldots,R_n)\Theta(R_1,\ldots,R_n)^*,
\end{equation}
where $P_\cE$ is the orthogonal projection of $F^2(H_n)\otimes \cK$ onto $\cE$.
According to Lemma \ref{inv.sub}, the subspace $\cN_J\otimes \cK$ is invariant
under the operators $R_i^*\otimes I_\cK$, $i=1,\ldots, n$.
Moreover, using the remarks preceding the theorem we have
$$\Theta(W_1,\ldots,W_n)=
P_{\cN_J\otimes \cK}\Theta(R_1,\ldots,R_n)|\cN_J\otimes \cK.$$
Hence, and compressing equation \eqref{proj}
to the subspace $\cN_J\otimes \cK$, we obtain
$$
P_{\cN_J\otimes \cK}P_\cE|\cN_J\otimes \cK=\Theta(W_1,\ldots,W_n)\Theta(W_1,\ldots,W_n)^*.
$$
Notice that, due to \eqref{E}, the left hand side of this equality is equal to  $P_\cM$, the orthogonal projection  of $\cN_J\otimes \cK$  onto $\cM$. Hence,
$$
P_\cM=\Theta(W_1,\ldots,W_n)\Theta(W_1,\ldots,W_n)^*
$$
and $\Theta(W_1,\ldots,W_n)$ is a partial isometry.
Therefore,
$$\cM=\Theta(W_1,\ldots, W_n)\left[ \cN_J\otimes \cG\right]$$
 and the proof
is complete.
\end{proof}

We remark that in the particular case when the ideal $J$ is generated by the polynomials
$S_iS_j-S_jS_i$, \ $i,j=1,\ldots, n$, then $\cN_J=F^2_s$, the symmetric Fock space,  and $B_i$, $i=1,\dots, n$, are the creation operators on the symmetric Fock space. In this case Theorem \ref{Beur} provides a Beurling type theorem for Arveson's space $H^2$  which was also obtained in \cite{McT} using different methods.

From now on, throughout this section,  we assume    that
$J$
is a WOT-closed two-sided ideal of $F_n^\infty$ such that $1\in\cN_J$.

\begin{theorem}\label{compact}
Let $J$
be a WOT-closed two-sided ideal of $F_n^\infty$ such that $1\in\cN_J$.
Then all the compact operators in $B(\cN_J)$ are contained in the operator
space
$$\overline{\text{\rm span}}\{B_\alpha B_\beta^*:\ \alpha,\beta\in \FF_n^+\}.
$$
Moreover,   the $C^*$-algebra  $C^*(B_1,\ldots, B_n)$ is irreducible.
\end{theorem}
\begin{proof}
Since $\cN_J$ is an  invariant subspace under each operator $S_i^*$, \ $i=1,\ldots,n$, and contains the constants,
we have
\begin{equation*}
\begin{split}
I_{\cN_J}-B_1B_1^*-\cdots -B_nB_n^*&=P_{\cN_J}(I-S_1S_1^*-\cdots -S_nS_n^*)|\cN_J\\
&=P_{\cN_J}P_\CC|\cN_J\\
&=P^{\cN_J}_\CC,
\end{split}
\end{equation*}
where $ P^{\cN_J}_\CC$ is the orthogonal projection of $\cN_J$ onto $\CC$.
Let
$$g(S_1,\ldots, S_n):=\sum\limits_{|\alpha|\leq m} a_\alpha S_\alpha\quad \text{ and } \quad \xi:=\sum\limits_{\beta\in \FF_n^+} b_\beta e_\beta\in \cN_J\subset F^2(H_n).
$$
Notice that
\begin{equation*}
\begin{split}
P^{\cN_J}_\CC g(B_1,\ldots, B_n)^*\xi&=\sum_{|\alpha|\leq m}P_\CC \overline{a}_\alpha S_\alpha^*\xi=\sum_{|\alpha|\leq m}\overline{a}_\alpha  b_\alpha
\\
&=\left< \xi, \sum_{|\alpha|\leq m} a_\alpha e_\alpha\right>=
\left< \xi,g(B_1,\ldots, B_n)(1)\right>.
 \end{split}
\end{equation*}
Therefore,
\begin{equation}\label{rankone}
f(B_1,\ldots, B_n)P^{\cN_J}_\CC g(B_1,\ldots, B_n)^*\xi=
\left< \xi,g(B_1,\ldots, B_n)(1)\right>f(B_1,\ldots, B_n)(1)
\end{equation}
for any polynomial $f(B_1,\ldots, B_n):=\sum\limits_{|\gamma|\leq p}c_\gamma B_\gamma$. Hence, $f(B_1,\ldots, B_n)P^{\cN_J}_\CC g(B_1,\ldots, B_n)^*$ is a rank one operator in $B(\cN_J)$.
Moreover, since
$P_\CC^{\cN_J}=I_{\cN_J}-B_1B_1^*-\cdots -B_nB_n^*$, we deduce that the above operator is also in the operator space $\overline{\text{\rm span}}\{B_\alpha B_\beta^*:\ \alpha,\beta\in \FF_n^+\}$.

Since  the polynomials $\sum\limits_{|\alpha|\leq m}a_\alpha e_\alpha$, \ $m\in \NN$, $a_\alpha\in \CC$,
are dense in $F^2(H_n)$, it is clear that the set
$$\cL:=\left\{\left(\sum_{|\alpha|\leq m}a_\alpha B_\alpha\right)(1):\ m\in \NN, a_\alpha\in \CC\right\}$$
is  dense in $\cN_J$.
Using this density and  relation \eqref{rankone}, we deduce that
all compact operators  in $B(\cN_J)$ are included in the operator space
$\overline{\text{\rm span}}\{B_\alpha B_\beta^*:\ \alpha,\beta\in \FF_n^+\}$.

To prove the last part of this lemma, let $\cM$ be a nonzero subspace of $\cN_J$ which is jointly  reducing
for $B_1,\ldots, B_n$. Take $f\in \cM$, $f\neq 0$ and assume that $f=a_0+\sum\limits_{|\alpha|\geq 1} a_\alpha e_\alpha.
$
If
 $a_\beta$  is a nonzero coefficient of $f$,
then
\begin{equation}\label{Bbeta}
B_\beta^*f=P_{\cN_J}S_\beta^*f=P_{\cN_J}\left(a_\beta+\sum_{|\gamma|\geq 1}a_{\beta\gamma }e_\gamma\right)
\end{equation}
 is  in $\cM$.
Since $1\in \cN_J$, we have $\left<P_{\cN_J}e_\gamma,1\right>=\left< e_\gamma, 1\right>=0$ for any $\gamma\in \FF_n^+$ with $|\gamma|\geq 1$. Hence, we deduce that
$$
P_{\cN_J}1=1\quad \text{ and }\quad P_\CC P_{\cN_J}\left(\sum_{|\gamma|\geq 1}a_{\beta\gamma }e_\gamma\right)=0.
$$
Therefore, relation \eqref{Bbeta} implies $P_\CC P_{\cN_J}f=a_\beta$.

On the other hand, since
$P_\CC^{\cN_J}=I_{\cN_J}-B_1B_1^*-\cdots -B_nB_n^*$  and $\cM$ is reducing for
$B_1,\ldots, B_n$, we infer that
$a_\beta\in \cM$, so $1\in \cM$. Using again that $\cM$ is invariant under $B_1,\ldots, B_n$, we have
 $\cL\subseteq \cM$. Since $\cL$ is dense in $\cN_J$, we deduce that $\cN_J\subset \cM$ and therefore $\cN_J= \cM$.
This completes the proof.
\end{proof}

We say that two row contractions $[T_1,\ldots, T_n]$, \ $T_i\in B(\cH)$, and $[T_1',\ldots, T_n']$, \ $T_i'\in B(\cH')$, are unitarily equivalent if there exists a unitary operator $U:\cH\to \cH'$ such that
$T_i=U^* T_i' U$ for  any $i=1,\ldots, n$.
If $[B_1,\ldots, B_n]$ is a constrained shift as above, then $[B_1\otimes I_\cH,\ldots, B_n\otimes I_\cH]$ is called constrained shift with multiplicity $\dim \cH$.

\begin{proposition}\label{eq-mult}
Two constrained shifts are unitarily equivalent if and only if their multiplicities are equal.
\end{proposition}
\begin{proof} Let $[B_1\otimes I_\cH,\ldots, B_n\otimes I_\cH]$ and $[B_1\otimes I_{\cH'},\ldots, B_n\otimes I_{\cH'}]$  be two
constrained shifts and let $U:\cN_J\otimes \cH\to \cN_J\otimes \cH'$ be a unitary  operator such that
$$
U(B_i\otimes I_\cH)=(B_i\otimes I_{\cH'})U,\quad i=1,\ldots, n.
$$
Since $U$ is unitary , we deduce that
$$
U(B_i^*\otimes I_\cH)=(B_i^*\otimes I_{\cH'})U,\quad i=1,\ldots, n.
$$
Since, according to Theorem \ref{compact}, the $C^*$-algebra $C^*(B_1,\ldots, B_n)$ is irreducible, we infer  that
$U=I_{\cN_J}\otimes W$ for some unitary operator $W\in B(\cH,\cH')$. Therefore, $\dim \cH=\dim \cH'$. The converse is obvious.
 \end{proof}

We need a few more definitions.
Let $\cS\subseteq B(\cK)$ be a set of operators acting on the Hilbert space $\cK$.
Denote by $A(\cS)$ the nonselfadjoint algebra generated by $S$ and the identity, and let $C^*(\cS)$ be the $C^*$-algebra generated by $\cS$ and the identity.
A  subspace $\cH\subseteq \cK$  is called
$*$-cyclic  for  $\cS$    if
$$\cK=\bigvee\{Xh:\ X\in C^*(\cS), \ h\in \cH\},
$$
 i.e., $\cK$ is the smallest reducing subpace for $\cS$ which contains $\cH$.
We call $\cH$   cyclic for $\cS$   if
$$\cK=\bigvee\{Xh:\ X\in A(\cS), \ h\in \cH\ \},
$$
 i.e., $\cK$ is the smallest invariant subpace under $\cS$ which contains $\cH$.
Finally, a  subspace $\cH\subseteq \cK$ is called  co-invariant under $\cS$ if
$X^*\cH\subseteq \cH$ for any $X\in \cS$.

\begin{theorem}\label{cyclic}
Let  $J$ be a WOT-closed two-sided ideal of $F_n^\infty$ such that $1\in \cN_J$ and let $\cD$  be a Hilbert space.
If $\cM\subseteq \cN_J\otimes \cD$ is a co-invariant subspace under $B_i\otimes I_\cD$, \ $i=1,\ldots, n$,
then  there exists a subspace $\cE\subseteq \cD$ such that
\begin{equation}\label{span}
\overline{\text{\rm span}}\,\left\{\left(B_\alpha\otimes I_\cD\right)\cM:\ \alpha\in \FF_n^+\right\}=\cN_J\otimes \cE.
\end{equation}
\end{theorem}

\begin{proof}
Denote  by $P_0:=P_\CC\otimes I_\cD$ the orthogonal projection of $\cN_J\otimes \cD$ onto  the subspace $1\otimes \cD$, which is identified with $\cD$. Define the subspace $\cE\subset \cD$ by setting $\cE:=P_0\cM$ and  let
$f$ be a nonzero element of $\cM$ having the Fourier representation
$$f=\sum_{\alpha\in \FF_n^+}e_\alpha\otimes h_\alpha,\quad h_\alpha\in \cD.
$$
Let $\beta\in \FF_n^+$ be such that $h_\beta\neq 0$ and note that
\begin{equation}\label{bb*}
\begin{split}
(B_\beta^*\otimes I_\cD)f&=
(P_{\cN_J}\otimes I_\cD)(S_\beta^*\otimes I_\cD)f\\
&=
(P_{\cN_J}\otimes I_\cD)\left( 1\otimes h_\beta+\sum_{|\gamma|\geq 1} e_\gamma\otimes h_{\beta\gamma}\right).
\end{split}
\end{equation}
As in the proof of Theorem \ref{compact}, since $1\in \cN_J$, we have $P_{\cN_J}1=1$ and $P_\CC P_{\cN_J}e_\gamma=0$ for $|\gamma|\geq 1$.
Hence and using  \eqref{bb*}, we obtain
$P_0(B_\beta^*\otimes I_\cD)f=h_\beta$. Since $\cM$ is  a co-invariant
subspace under $B_i\otimes I_\cD$, \ $i=1,\ldots, n$,  it is clear that
$h_\beta\in \cE$.  Using this and  taking into account that $1\in\cN_J$, we  deduce that
$$
(B_\beta\otimes I_\cD)(1\otimes h_\beta)=P_{\cN_J}e_\beta\otimes h_\beta\in \cN_J\otimes \cE.
$$
Now, since $f\in \cM\subseteq \cN_J\otimes \cD$, we infer that
$$
f=(P_{\cN_J}\otimes I_\cD)f=\lim_{k\to\infty}\sum_{|\alpha|=k} P_{\cN_J}e_\beta\otimes h_\beta
$$
is in $\cN_J\otimes \cE$. This shows that $\cM\subset\cN_J\otimes \cE$ and
therefore
$$
\cY:= \overline{\text{\rm span}}\,\left\{\left(B_\alpha\otimes I_\cD\right)\cM:\ \alpha\in \FF_n^+\right\}\subset\cN_J\otimes \cE.
$$

For the other inclusion, we prove first that $\cE\subset \cY$. If $h_0\in \cE$, $h_0\neq 0$, then there exists $g\in \cM$ such that
$g=1\otimes h_0+\sum\limits_{|\alpha|\geq 1}  e_\alpha\otimes h_\alpha.$
Due to the first part of the proof, we have
$$
h_0=P_0 g=\left(I-\sum_{i=1}^n (B_i\otimes I_\cD) (B_i\otimes I_\cD)^*\right)g.
$$
Acoording to the proof of  Theorem \ref{compact}, we have $P_\CC^{\cN_J}=I_{\cN_J}-\sum_{i=1}B_iB_i^*$.
Using this and the fact that $\cM$ is co-invariant under
$B_i\otimes I_\cD$, \ $i=1,\ldots, n$,  we deduce that $h_0\in\cY$ for any $h_0\in \cE$, i.e., $\cE\subset \cY$. The latter inclusion  shows that
$(B_\alpha\otimes I_\cD)(1\otimes \cE)\subset \cY$ for any
$\alpha\in \FF_n^+$, which implies
\begin{equation}\label{PN}
P_{\cN_J}e_\alpha\otimes \cE\subset \cY,\quad \alpha\in \FF_n^+.
\end{equation}
Let
$\varphi\in \cN_J\otimes \cE\subset F^2(H_n)\otimes \cE$ be with Fourier representation
$\varphi=\sum\limits_{\alpha\in \FF_n^+}e_\alpha\otimes k_\alpha,\quad k_\alpha\in \cE.
$
Using   relation \eqref{PN}, we have
$$
\varphi=(P_{\cN_J}\otimes I_\cE)g=\lim_{k\to \infty}\sum_{|\alpha|=k}P_{\cN_J}e_\alpha\otimes k_\alpha\in \cY.
$$
Therefore, $\cN_J\otimes \cE\subseteq \cY$, which
 completes the proof.
\end{proof}

\begin{corollary}\label{cyc} Let  $J$ be a WOT-closed two-sided ideal of $F_n^\infty$ such that $1\in \cN_J$ and let $\cD$  be a Hilbert space.
Let $\cM\subseteq \cN_J\otimes \cD$  be a co-invariant subspace  under $B_i\otimes I_\cD$, \ $i=1,\ldots, n$.
 Then the following statements  are equivalent:
\begin{enumerate}
\item[(i)] $\cM$ is a cyclic subspace for $B_i\otimes I_\cD$, \ $i=1,\ldots, n$;
\item[(ii)] $P_0 \cM=\cD$;
\item[(iii)]$\cM^\perp\bigcap \cD=\{0\}$.
\end{enumerate}
\end{corollary}

\begin{proof}
The equivalence $(i)\leftrightarrow (ii)$ is clear from Theorem \ref{cyclic} and the definition of cyclic subspace.
To prove that $(ii)\leftrightarrow (iii)$, notice first  that if there exists
$h\in \cM^\perp\bigcap \cD$, \ $h\neq 0$, then taking into account that $\cM^\perp$ is invariant under
$B_i\otimes I_\cD$, \ $i=1,\ldots, n$, we deduce that $\cN_J\otimes h\subset \cM^\perp$. This shows that $h\notin P_0\cM$ which means that
$P_0 \cM$ is not equal to $\cD$. Now, assume that there exists $k\in \cD$, $k\neq 0$, such that $k\perp P_0\cM$. Since $1\in \cN_J$, we have $k\perp \cM$ which shows that $k\in \cD\bigcap \cM^\perp$. The proof is complete.
\end{proof}

\begin{corollary}
Let  $J$ be a WOT-closed two-sided ideal of $F_n^\infty$ such that $1\in \cN_J$ and let $\cD$  be a Hilbert space.
 A subspace $\cM\subseteq \cN_J\otimes \cD$ is reducing under each operator $B_i\otimes I_\cD$, \ $i=1,\ldots, n$,
if and only if   there exists a subspace $\cE\subseteq \cD$ such that
\begin{equation*}
 \cM=\cN_J\otimes \cE.
\end{equation*}
\end{corollary}

The next result, is a Wold type decomposition for nondegenerate $*$-representations of the $C^*$-algebra
$C^*(B_1,\ldots, B_n)$.

\begin{theorem}\label{wold}
Let  $J$ be a WOT-closed two-sided ideal of $F_n^\infty$ such that $1\in \cN_J$, and
let $\pi:C^*(B_1,\ldots, B_n)\to B(\cK)$ be a nondegenerate $*$-representation  of $C^*(B_1,\ldots, B_n)$ on a separable  Hilbert space  $\cK$. Then
$\pi$ decomposes into
a direct sum
$$
\pi=\pi_0\oplus \pi_1 \  \text{ on  } \ \cK=\cK_0\oplus \cK_1,
$$
where $\pi_0$, $\pi_1$  are disjoint representations of $C^*(B_1,\ldots, B_n)$
on the Hilbert spaces
$$
\cK_0:=\overline{\text{\rm span}}\left\{\pi(B_\alpha)
\left(I-\sum_{i=1}^n \pi(B_i)\pi(B_i)^*\right)
\cK:\  \alpha\in \FF_n^+\right\}\quad  \text{ and } \quad \cK_1:=\cK_0^\perp,
$$
 respectively, such that, up to an isomorphism,
\begin{equation}\label{shi}
\cK_0\simeq\cN_J\otimes \cG, \quad  \pi_0(X)=X\otimes I_\cG, \quad  X\in C^*(B_1,\ldots, B_n),
\end{equation}
 for some Hilbert space $\cG$ with
$$
\dim \cG=\dim \left[\text{\rm range}\,\left(I-\sum_{i=1}^n \pi(B_i)\pi(B_i)^*\right)\right],
$$
 and $\pi_1$ is a $*$-representation  which annihilates the compact operators   and
$$
\pi_1(B_1)\pi_1(B_1)^*+\cdots +\pi_1(B_n)\pi_1(B_n)^*=I_{\cK_1}.
$$
Moreover, if $\pi'$ is another nondegenerate  $*$-representation of
$C^*(B_1,\ldots, B_n)$ on a separable  Hilbert space  $\cK'$, then $\pi$ is unitarily equivalent to $\pi'$ if and only if $\dim\cG=\dim\cG'$ and $\pi_1$ is unitarily equivalent to $\pi_1'$.
\end{theorem}

\begin{proof}
Since   the subspace $\cN_J$ contains the constants, Theorem
\ref{compact} implies that
all the compact operators $\cL\cC(\cN_J))$ in $B(\cN_J)$ are contained in
$C^*(B_1,\ldots, B_n)$.
According to the standard theory of representations of the $C^*$-algebras, the representation
$\pi$ decomposes into
a direct sum
$\pi=\pi_0\oplus \pi_1$ on  $ \cK=\cK_0\oplus \cK_1$,
where
$$\cK_0:=\overline{\text{\rm span}}\{\pi(X)\cK:\ X\in \cL\cC(\cN_J)\}\quad \text{ and  }\quad  \cK_1:=\cK_0^\perp,
$$
and the the representations $\pi_j:C^*(B_1,\ldots, B_n)\to \cK_j$ are defined by $\pi_j(X):=\pi(X)|\cK_j$, \ $j=0,1$.
Now, it is clear that
$\pi_0$, $\pi_1$  are disjoint representations of $C^*(B_1,\ldots, B_n)$
 such that
 $\pi_1$ annihilates  the compact operators in $B(\cN_J)$, and  $\pi_0$ is uniquely determined by the action of $\pi$ on the ideal $\cL\cC(\cN_J))$.
Since every representation of $\cL\cC(\cN_J))$ is equivalent to a multiple of the identity representation
(see \cite{Arv-book}), we deduce
\eqref{shi}.
Now, we show that the space $\cK_0$ coincides with the one defined in the theorem. Using Theorem \ref{compact} and its proof, we deduce that
\begin{equation*}\begin{split}
\cK_0&:=\overline{\text{\rm span}}\{\pi(X)\cK:\ X\in \cL\cC(\cN_J)\}\\
&=\overline{\text{\rm span}}\{\pi(B_\alpha P_\CC^{\cN_J} B_\beta^*)\cK:\ \alpha, \beta\in \FF_n^+\}\\
&=
\overline{\text{\rm span}}\left\{\pi(B_\alpha)
\left(I-\sum_{i=1}^n \pi(B_i)\pi(B_i)^*\right)
\cK:\  \alpha\in \FF_n^+\right\}.
\end{split}
\end{equation*}

On the other hand, since $P_\CC^{\cN_J}=I-\sum\limits_{i=1}^n B_iB_i^*$ is a rank one projection in
$C^*(B_1,\ldots, B_n)$ (see Theorem \ref{compact}), we deduce that
$\sum\limits_{i=1}^n \pi_1(B_i)\pi_1(B_i)^*=I_{\cK_1}$, and
$$
\dim \cG=\dim \left[\text{\rm range}\,\pi(P_\CC^{\cN_J})\right]
=\dim \left[\text{\rm range}\,\left(I-\sum_{i=1}^n \pi(B_i)\pi(B_i)^*\right)\right].
$$

To prove the uniqueness, note that
according to the standard theory of representations of  $C^*$-algebras,
$\pi$ and $\pi'$ are unitarily equivalent if and only if
 $\pi_0$ and $\pi_0'$ (resp.~$\pi_1$ and $\pi_1'$) are unitarily equivalent. Using Proposition \ref{eq-mult}, we deduce that $\dim\cG=\dim\cG'$ and
 complete the proof.
\end{proof}

\begin{corollary}\label{Wold}
Under the hypotheses and notations of Theorem $\ref{wold}$ and setting
$V_i:=\pi(B_i)$, \ $i=1,\ldots, n$,
we have:
\begin{enumerate}
\item[(i)] $Q:=I_\cK-\sum\limits_{i=1}^n V_iV_i^*$ is an orthogonal projection and
$Q\cK=\bigcap\limits_{i=1}^n \ker V_i^*$;
\item[(ii)] $\cK_0=\left\{ k\in \cK:\ \lim\limits_{k\to \infty}\sum\limits_{|\alpha|=k}\|V_\alpha^*k\|^2=0\right\}$;
\item[(iii)]$\cK_1=\left\{ k\in \cK:\ \sum\limits_{|\alpha|=k}\|V_\alpha^*k\|^2=\|k\|^2 \ \text{ for any } k=1,2, \ldots\right\}$;
\item[(iv)] $\text{\rm SOT-}\lim\limits_{k\to\infty}\sum\limits_{|\alpha|=k}V_\alpha V_\alpha^*=P_{\cK_1}$;
\item[(v)]
$\sum\limits_{k=0}^\infty\sum\limits_{|\alpha|=k}V_\alpha QV_\alpha^*=P_{\cK_0}$;
\end{enumerate}
\end{corollary}

\begin{proof}
Since $I_{\cN_J}-\sum\limits_{i=1}^n B_iB_i^*=P_\CC^{\cN_J}$ is an orthogonal projection (see the proof of Theorem \ref{compact}), so is
$Q=\pi(P_\CC^{\cN_J})$. Therefore,
\begin{equation*}
\begin{split}
Q\cK&=\left\{ k\in \cK:\ \left(I-\sum_{i=1}^n V_iV_i^*\right)k=k\right\}\\
&=\left\{ k\in \cK:\ \sum_{i=1}^n V_iV_i^*k=0\right\}\\
&=\bigcap_{i=1}^n \ker V_i^*,
\end{split}
\end{equation*}
which proves (i).
Using  Theorem \ref{wold}, we have
\begin{equation}\label{matr}
\sum\limits_{|\alpha|=k} V_\alpha V_\alpha^*=
\left[
\begin{matrix}
\sum\limits_{|\alpha|=k} B_\alpha B_\alpha^*\otimes I_\cG &0\\
0&I_{\cK_1}
\end{matrix}\right].
\end{equation}
  Since
$\cN_J$ is co-invariant under $S_1,\ldots, S_n$, and
$B_i^*=S_i^*|\cN_J$, $i=1,\ldots, n$, we have
$$
\text{\rm SOT-}\lim_{k\to \infty}\sum\limits_{|\alpha|=k} B_\alpha B_\alpha^*\otimes I_\cG=
\text{\rm SOT-}\lim_{k\to \infty}
\left[P_{\cN_J}\left(\sum\limits_{|\alpha|=k} S_\alpha S_\alpha^*\right)|\cN_J\right]\otimes I_\cG=0.
$$
The latter equality holds due to the fact that
$[S_1,\ldots, S_n]$ is a pure contraction. Therefore, (iv) holds. Hence, and taking into account that
$$\sum_{i=1}^m \sum\limits_{|\alpha|=k} V_\alpha  QV_\alpha^*=I-\sum\limits_{|\alpha|=m+1} V_\alpha V_\alpha^*,
$$
we deduce (v).
Now, let $k\in \cK=\cK_0\oplus \cK_1$,\ $k=k_0+ k_1$, with
$k_0\in \cK_0$ and $k_1\in \cK_1$. By \eqref{matr}, we have
 \begin{equation}
\label{VVBB}
\sum\limits_{|\alpha|=m} \|V_\alpha^*k\|^2=\left<\left(
\sum\limits_{|\alpha|=m} B_\alpha B_\alpha^*\otimes I_\cG\right)k_0, k_0\right>+\|k_1\|^2, \quad m=1,2,\ldots.
\end{equation}
Hence, $\lim\limits_{m\to \infty}\sum\limits_{|\alpha|=m}\|V_\alpha^*k\|^2=0 $ if and only if
$k_1=0$, i.e., $k=k_0\in \cK_0$, which proves (ii).
On the other hand, \eqref{VVBB} shows that
$\sum\limits_{|\alpha|=m}\|V_\alpha^*k\|^2=\|k\|^2$  for any $m=1,2,\ldots$, if and only if
$\left<\left(
\sum\limits_{|\alpha|=m} B_\alpha B_\alpha^*\otimes I_\cG\right)k_0, k_0\right>=\|k_0\|^2$
 for any $m=1,2,\ldots$. Since
$[B_1,\ldots, B_n]$ is a pure row contraction,
$\text{\rm SOT-}\lim\limits_{m\to\infty} \sum\limits_{|\alpha|=m} B_\alpha B_\alpha^*=0$. Therefore, the above  equality holds for any $m=1,2,\ldots$, if and only if $k_0=0$, which is equivalent to $k=k_1\in \cK_1$. This completes the proof.
\end{proof}

\begin{corollary}\label{co-shi}
Let $\pi$  be a nondegenerate
$*$-representation of $C^*(B_1,\ldots, B_n)$ on a separable Hilbert space $\cK$, and let $V_i:=\pi(B_i)$, \ $i=1,\ldots, n$.
Then  the following statements are equivalent:
\begin{enumerate}
\item[(i)]
$V:=[V_1,\ldots, V_n]$ is a constrained shift;
 \item[(ii)]
$
\cK=\overline{\text{\rm span}}\left\{V_\alpha
\left(I-\sum\limits_{i=1}^n V_i V_i^*\right)
\cK:\  \alpha\in \FF_n^+\right\};
$
\item[(iii)] \text{\rm SOT-}$\lim\limits_{k\to\infty}\sum\limits_{|\alpha|=k} V_\alpha V_\alpha^*=0$.
\end{enumerate}
In this case,   the multiplicity of $V$ $($denoted by $\text{\rm mult}\,(V))$
 satisfies the equality
\begin{equation}\label{multiplicity}
\text{\rm mult}\,(V)=\dim (I-V_1V_1^*-\cdots -V_nV_n^*) \cK,
\end{equation}
and  it is  also equal to the minimum dimension of a cyclic subspace for $V_1,\ldots, V_n$.
\end{corollary}

\begin{proof}
The  above equivalences are consequences of Corollary \ref{Wold}, and  relation \eqref{multiplicity} follows from Theorem
\ref{wold}. We prove the last part of the corollary.
According to (ii) and Corollary \ref{Wold}, $$
\cL:=\bigcap_{i=1}^n V_i^*=(I-V_1V_1^*-\cdots -V_nV_n^*) \cK
$$
 is a cyclic subspace for $V_1,\ldots, V_n$. Now, let $\cE$ be any cyclic subspace for $V_1,\ldots, V_n$, i.e., $\cK=\bigvee\limits_{\alpha\in \FF_n^+} V_\alpha \cE$, and denote $A:=P_\cL|\cE\in B(\cE, \cL)$, where $P_\cL$ is the orthogonal projection of $\cK$ onto $\cL$.
Assume that $k\in \cL\ominus T\cE$ and let $h\in \cE$.
Notice that
$$
\left<h,k\right>=\left<h,P_\cL k\right>=\left<Ah,k\right>=0.
$$
On the other hand, since $V_\alpha^*k=0$ for any $k\in \cL$, we have $\left<V_\alpha h, k\right>=0$ for any $\alpha\in \FF_n^+$ with $|\alpha|\geq 1$.
Therefore,  $k\perp V_\alpha\cE$ for any $\alpha\in \FF_n^+$. Since $\cE$ is a cyclic  subspacefor $V_1,\ldots, V_n$, we deduce that $k=0$, and therefore $\overline {A\cE}=\cL$. This shows that $A^*\in B(\cL, \cE)$ is one-to-one
and, consequently, $\dim\cL\leq \dim\cE$. This completes the proof.
\end{proof}

An easy consequence of Corollary \ref{co-shi} and Proposition \ref{eq-mult} is the following.
 \begin{proposition}
Two constrained shifts are similar if and only if  they are unitarily equivalent.
\end{proposition}
\begin{proof} One implication is obvious.
Let $V:=[V_1,\ldots, V_n]$, $V_i\in B(\cK)$, and
$V':=[V_1',\ldots, V_n']$, $V_i'\in B(\cK')$, be two constrained shifts and let $X:\cK\to \cK'$ be an invertible operator such that
$$XV_i=V_i'X,\quad i=1,\ldots, n.
$$
If $\cM$ is a cyclic subspaces for $V_1,\ldots, V_n$, then
\begin{equation*}
\begin{split}
\cK'&=X\cK=X\left(\bigvee_{\alpha\in \FF_n^+} V_\alpha \cM\right)\\
&\subseteq \bigvee_{\alpha\in \FF_n^+}X V_\alpha \cM=
\bigvee_{\alpha\in \FF_n^+} V_\alpha'X \cM\subseteq \cK'.
\end{split}
\end{equation*}
Therefore, $\cK'=\bigvee_{\alpha\in \FF_n^+} V_\alpha'X \cM$, which shows that $X\cM$ is cyclic for $V'$. Since $X$ is invertible, $\dim \cM=\dim X\cM$. Hence an using  Corollary \ref{co-shi}, we conclude that  the two constrained shifts have the same multiplicity.
By Proposition \ref{eq-mult}, the result follows.
\end{proof}

  \bigskip

      \section{Dilations for  constrained row contractions}\label{Dilations}

In this section we develop a dilation theory for row contractions $T:=[T_1,\ldots, T_n]$
subject to constraints such as
$$
p(T_1,\ldots, T_n)=0,\quad p\in \cP,
$$
where $\cP$ is a set of noncommutative polynomials. The model $n$-tuple is the universal row contraction $[B_1,\ldots, B_n]$ satisfying the same constraints as $T$.
 The theory is based on a class of
{\it constrained Poisson kernels} associated with $T$  and representations of the $C^*$-algebra generated  by $B_1,\ldots, B_n$ and the identity.
Under natural conditions  on the constraints we have uniqueness for the minimal dilation. We introduce  and evaluate
the {\it dilation index}, a numerical invariant for row contractions,  and show that it does not depend on the constraints.
These results are used  in the next sections in connection with
characteristic  functions and models for constrained  row contractions.

We need to recall from \cite{Po-poisson} a few facts about noncommutative
Poisson transforms associated with row contractions $T:=[T_1,\ldots, T_n]$,
\ $T_i\in B(\cH)$. For  each $0<r\leq 1$, define the defect operator
$\Delta_{T,r}:=(I-r^2T_1T_1^*-\cdots -r^2 T_nT_n^*)^{1/2}$.
The Poisson  kernel associated with $T$ is the family of operators
$$
K_{T,r} :\cH\to F^2(H_n)\otimes \overline{\Delta_{T,r}\cH}, \quad  0<r\leq 1,
$$
defined by

\begin{equation}\label{Poiss}
K_{T,r}h:= \sum_{k=0}^\infty \sum_{|\alpha|=k} e_\alpha\otimes r^{|\alpha|}
\Delta_{T,r} T_\alpha^*h,\quad h\in \cH.
\end{equation}
When $r=1$, we denote $\Delta_T:=\Delta_{T,1}$ and $K_T:=K_{T,1}$.
The operators $K_{T,r}$ are isometries if $0<r<1$, and
$$
K_T^*K_T=I-
\text{\rm SOT-}\lim_{k\to\infty} \sum_{|\alpha|=k} T_\alpha T_\alpha^*.
$$
This shows that $K_T$ is an isometry if and only if $T$ is a {\it pure} row
 contraction (\cite{Po-isometric}),
i.e.,
$$
\text{\rm SOT-}\lim_{k\to\infty} \sum_{|\alpha|=k} T_\alpha T_\alpha^*=0.
$$
A key property of the Poisson kernel
is that
\begin{equation}\label{eq-ker}
K_{T,r}(r^{|\alpha|} T_\alpha^*)=(S_\alpha^*\otimes I)K_{T,r}\qquad \text{ for any } 0<r\leq 1,\ \alpha\in \FF_n^+.
\end{equation}
  In \cite{Po-poisson}, we introduced  the Poisson transform associated with
 $T:=[T_1,\ldots, T_n]$ as the unital completely contractive  linear map $\Psi_T:C^*(S_1,\ldots, S_n)\to B(\cH)$ defined by
 \begin{equation*}
 \Psi_T(f):=\lim_{r\to 1} K_{T,r}^* (f\otimes I)K_{T,r},
\end{equation*}
 where the limit exists in the norm topology of $B(\cH)$. Moreover, we have
 $$
 \Psi_T(S_\alpha S_\beta^*)=T_\alpha T_\beta^*, \quad \alpha,\beta\in \FF_n^+.
 $$

 When $T$ is a completely non-coisometric (c.n.c.) row-contraction, i.e.,
 there is no $h\in \cH$, $h\neq 0$, such that
 $$
 \sum_{|\alpha|=k}\|T_\alpha^* h\|^2=\|h\|^2
 \quad \text{\rm for any } \ k=1,2,\ldots,
 $$
an  $F_n^\infty$-functional calculus was developed  in \cite{Po-funct}.
 We showed that if $f=\sum\limits_{\alpha\in \FF_n^+} a_\alpha S_\alpha$ is
 in $F_n^\infty$, then
 \begin{equation*}
 \Gamma_T(f)=f(T_1,\ldots, T_n):=
 \text{\rm SOT-}\lim_{r\to 1}\sum_{k=0}^\infty
  \sum_{|\alpha|=k} r^{|\alpha|} a_\alpha T_\alpha
\end{equation*}
exists and $\Gamma_T:F_n^\infty\to B(\cH)$ is a WOT-continuous
completely contractive homomorphism.
More about   noncommutative Poisson transforms  on $C^*$-algebras generated
by isometries can be found in
\cite{Po-poisson}, \cite{APo}, \cite{Po-tensor}, \cite{Po-curvature}, and \cite{Po-similarity}.

Throughout this section we keep the notations from Section \ref{Beurling}.
Let $J\neq F_n^\infty$ be a WOT-closed two-sided ideal of $F_n^\infty$ generated by a family
of polynomials $\cP_J\subset F_n^\infty$ and let \linebreak $T:=[T_1,\ldots, T_n]$, \ $T_i\in B(\cH)$, be  a row contraction such that
$$
p(T_1,\ldots, T_n)=0\quad \text{ for any }\quad p\in \cP_J.
$$
Let $\cD$ and $\cK$  be Hilbert spaces  and let
 $Z_i\in B (\cK)$ be  bounded operators such that
$[Z_1,\ldots, Z_n]$ is a Cuntz row contraction, i.e.,
$$Z_1Z_1^*+\cdots +Z_nZ_n^*=I_\cK. $$
An $n$-tuple $V:=[V_1,\ldots, V_n]$ of operators  with
\begin{equation}\label{Vi}
V_i:=\left[\begin{matrix}
B_i\otimes I_\cD&0\\
0&Z_i
\end{matrix}
\right], \qquad i=1,\ldots, n,
\end{equation}
where the $n$-tuple $[B_1,\ldots, B_n]$ is the  constrained shift associated with $J$,
is called  {\it constrained} (or $J$-{\it constrained}) {\it dilation} of $T$ if:
\begin{enumerate}
\item[(i)] $p(V_1,\ldots, V_n)=0$ for any $p\in \cP_J$;
\item[(ii)] $\cH$ can be identified with a co-invariant subspace under  $V_1,\ldots, V_n$ such that
$$
T_i=P_\cH V_i|\cH,\quad i=1,\ldots, n.
$$
\end{enumerate}
The dilation is minimal if $\cH$ is cyclic for $V_1,\ldots, V_n$, i.e.,
$$
(\cN_J\otimes \cD)\oplus \cK=\bigvee_{\alpha\in \FF_n^+} V_\alpha \cH.
$$
We introduce the {\it dilation index} of $T$,  denoted by $\text{\rm dil-ind}\,(T)$, to be the minimum dimension of $\cD$ such that $V$ is  a constrained   dilation of $T$.

Our first dilation result  for constrained row contractions is the following.

\begin{theorem}\label{dil1}
Let $J\neq F_n^\infty$ be a WOT-closed two-sided ideal of $F_n^\infty$ generated by a family
of polynomials $\cP_J\subset F_n^\infty$ and let $T:=[T_1,\ldots, T_n]$, \ $T_i\in B(\cH)$, be a row contraction such that
\begin{equation}\label{constraints1}
p(T_1,\ldots, T_n)=0\quad \text{ for any }\quad p\in \cP_J.
\end{equation}
Then there exists a Hilbert space $\cK$  and some operators
$Z_i\in B (\cK)$  with the properties
$$Z_1Z_1^*+\cdots +Z_nZ_n^*=I_\cK $$  and
$$
p(Z_1,\ldots, Z_n)=0\quad \text{ for any }  \quad p\in \cP_J,
$$
such that:
\begin{enumerate}
\item[(i)]
 $\cH$ can be identified with a  co-invariant subspace of  $\tilde\cK:=(\cN_J\otimes \overline{\Delta_T\cH})\oplus \cK$ under the operators
$$
V_i:=\left[\begin{matrix}
B_i\otimes I_{\overline{\Delta_T\cH}}&0\\0&Z_i
\end{matrix}\right],\quad i=1,\ldots,n;
$$
\item[(ii)]
$T_i^*=V_i^*|\cH$,\ $i=1,\ldots,n$.
\end{enumerate}
Moreover,
$\cK=\{0\}$ if and only if $[T_1,\ldots, T_n]$ is a pure row contraction.
\end{theorem}

\begin{proof}
Consider the subspace
$$
\cM:=\overline{\text{\rm span}}\left\{S_\alpha p(S_1,\ldots, S_n) S_\beta(1):\ p\in \cP_J, \alpha,\beta\in \FF_n^+\right\}.
$$
It is clear that $\cM\subseteq \cM_J$. To prove that  $\cM_J\subseteq \cM$, it is enough to show that $\cM^\perp\subseteq \cM_J^\perp$.
Let $g\in F^2(H_n)$  be such that
$$\left< g, S_\alpha p(S_1,\ldots, S_n) S_\beta(1)\right>=0\quad \text{ for any }\ p\in \cP_J, \alpha,\beta\in \FF_n^+.
$$
It is known (see \cite{APo}, \cite{DP}) that for   any $\varphi(S_1,\ldots, S_n)\in F_n^\infty$, there is a sequence of polynomials $\{q_m(S_1,\ldots, S_n)\}_{m=1}^\infty$ which is SOT-convergent
to  $\varphi(S_1,\ldots, S_n)$ as $m\to\infty$.
Consequently,
$$\left< g,\varphi(S_1,\ldots, S_n) p(S_1,\ldots, S_n) S_\beta(1)\right>=0
$$ for any  $\varphi(S_1,\ldots, S_n)\in F_n^\infty$, $p\in \cP_J$, and $\alpha,\beta\in \FF_n^+.$
Hence, $g\in \cM_J^\perp$. Therefore, $\cM_J=\cM$.

Now, using the properties of the Poisson kernel
$K_T$ (see \eqref{eq-ker}) and that
$p(T_1,\ldots, T_n)=0$  for any $p\in \cP_J$, we obtain
$$
\left<  K_Tk, S_\alpha p(S_1,\ldots, S_n) S_\beta(1)\otimes h\right>=\left<k,T_\alpha p(T_1,\ldots,T_n)T_\beta \Delta_Th\right>=0
$$
for any $k\in \cH$, $h\in \overline{\Delta_T\cH}$, and $p\in \cP_J$.
Since $\cM_J=\cM$, we infer that
\begin{equation}\label{range}
\text{\rm range}\,K_T\subseteq\left(\cM_J\otimes \overline{\Delta_T\cH}
\right)^\perp =\cN_J\otimes \overline{\Delta_T\cH}.
\end{equation}

Consider the {\it constrained Poisson kernel} $K_{J,T}:\cH\to \cN_J\otimes \overline{\Delta_T \cH}$ defined by
$$K_{J,T}:=(P_{\cN_J}\otimes I_{\overline{\Delta_T \cH}})K_T,
$$
where $K_T$ is the Poisson kernel defined by \eqref{Poiss}.
Using relations \eqref{eq-ker} and \eqref{range}, we obtain
\begin{equation}\label{KJ}
K_{J,T}T_\alpha^*=(B_\alpha^*\otimes I_\cH)K_{J,T},\quad  \alpha\in \FF_n^+.
\end{equation}
Define the contraction $Q:=\text{\rm SOT-}\lim\limits_{k\to\infty} \sum\limits_{|\alpha|=k}
T_\alpha T_\alpha^*$ and the operator
$$Y:\cH\to \cK:=\overline{Q^{1/2}\cH}\quad \text{ by }\quad
Yh:=Q^{1/2}h, \ h\in \cH.
$$
For each $i=1,\ldots, n$, define $\Lambda_i:Q^{1/2}\cH\to \cK$ by  setting
\begin{equation}\label{Z_i}
\Lambda_i Yh:=YT_i^*h,\quad h\in \cH.
\end{equation}
The operators $\Lambda_i$, $i=1,\ldots, n$, are well-defined since
$$
\sum_{i=1}^n\|\Lambda_i Yh\|^2=\left< \sum_{i=1}^n T_i QT_i^*h,h\right>=\|Q^{1/2}h\|^2.
$$
Therefore the operator $\Lambda_i$ can be extended  to a bounded operator on
$\cK$, which will also be denoted by $\Lambda_i$.
Now,  setting $Z_i:=\Lambda_i^*$, $i=1,\ldots, n$, relation  \eqref{Z_i} implies
\begin{equation}
\label{YZT}
Y^*Z_i=T_iY^*,\quad i=1,\ldots, n.
\end{equation}
Notice
that
\begin{equation*}
\begin{split}
Y^*\left( \sum_{i=1}^n Z_iZ_i^*\right) Y&= \sum_{i=1}^nT_i Y^*Y T_i^*\\
&=\sum_{i=1}^nT_i Q  T_i^*=Q=YY^*.
\end{split}
\end{equation*}
Hence,
$$\left<\left(\sum_{i=1}^n Z_iZ_i^*\right)Yh,Yh\right>=
\left< Yh,Yh\right>,
$$
which implies $\sum\limits_{i=1}^n Z_iZ_i^*=I_\cK$. Using relations \eqref{YZT} and \eqref{constraints1}, we get
$$Y^*p(Z_1,\ldots, Z_n)=p(T_1,\ldots, T_n)Y^*=0, \quad p\in\cP_J.
$$
Since $Y^*$ is injective on $\cK=\overline{Y\cH}$, we deduce that
$p(Z_1,\ldots, Z_n)=0$ for any $p\in\cP_J$.

Define the operator $V:\cH\to [\cN_J\otimes \cH]\oplus \cK$ by setting
$$V:=\left[\begin{matrix}
K_{J,T}\\ Y
\end{matrix}\right].
$$
Note that
\begin{equation*}\begin{split}
\|Vh\|^2&=\|K_{J,Y}h\|^2+\|Yh\|^2\\
&\|h\|^2-\text{\rm SOT-}\lim_{k\to\infty}\left<\sum_{|\alpha|=k}T_\alpha T_\alpha^*\right>+\|Yh\|^2\\
&=\|h\|^2-\left<Qh,h\right> + \left<Qh,h\right>\\
&=\|h\|^2
\end{split}
\end{equation*}
for any $h\in \cH$. Therefore, $V$ is an isometry.
On the other hand, using relations \eqref{KJ} and \eqref{Z_i}, we deduce that
\begin{equation*}
\begin{split}
VT_i^*&=\left[\begin{matrix}
K_{J,T}\\ Y
\end{matrix}\right]T_i^*h=K_{J,T}T_i^*h\oplus YT_i^*h\\
&=(B_i^*\otimes I_\cH)K_{J,T}h\oplus Z_i^*Yh\\
&=\left[\begin{matrix}
B_i^*\otimes I_{\overline{\Delta_T\cH}}&0\\0&Z_i^*
\end{matrix}\right]Vh.
\end{split}
\end{equation*}
Since $V$ is an isometry we can identify  $\cH$  with $V\cH$ and complete the proof of (i) and (ii). The  last part  of the theorem is obvious.
\end{proof}

\begin{corollary}
In the particular case when $n=1$ and $\cP_J=0$, we obtain
the classical isometric dilation theorem for contractions
obtained by Sz.-Nagy $($see \cite{SzF-book}$)$ by different methods.
\end{corollary}

Now  we can evaluate the dilation index of a constrained row contraction and show that it does not depend on
the constraints. Moreover, we show that the dilation index coincides with $\rank \Delta_T$, which is also called the pure rank of $T$ in \cite{DKS}.

\begin{corollary}\label{dil-ind}
Let $J\neq F_n^\infty$ be a WOT-closed two-sided ideal of $F_n^\infty$ generated by a family
of polynomials $\cP_J\subset F_n^\infty$ and let $T:=[T_1,\ldots, T_n]$, \ $T_i\in B(\cH)$, be a row contraction such that
\begin{equation*}
p(T_1,\ldots, T_n)=0\quad \text{ for any }\quad p\in \cP_J.
\end{equation*}
Then the dilation index of $T$  satisfies the equation
$$
\text{\rm dil-ind}\,(T)=\rank \Delta_T.
$$
\end{corollary}

\begin{proof}
Let $\cD$ and $\cK$  be Hilbert spaces  and let
 $Z_i\in B (\cK)$ be  bounded operators such that
$$Z_1Z_1^*+\cdots +Z_nZ_n^*=I_\cK. $$
Assume that   the $n$-tuple $V:=[V_1,\ldots, V_n]$ given by
\begin{equation}\label{ViBi}
V_i:=\left[\begin{matrix}
B_i\otimes I_\cD&0\\
0&Z_i
\end{matrix}
\right], \qquad i=1,\ldots, n,
\end{equation}
is a  constrained dilation of $T$. Since $\cH$ is co-invariant under $V_1,\ldots, V_n$, and $\cN_J$ is co-invariant under the left creation operators $S_1,\ldots, S_n$,  we have
$$
I_\cH-\sum_{i=1}^n T_iT_i^*=P_\cH\left[
\begin{matrix}
\left[P_{\cN_J}\left( I-\sum\limits_{i=1}^n S_iS_i^*\right)|\cN_J\right]\otimes I_\cD& 0\\
0&0
\end{matrix}\right]|\cH.
$$
Hence,  and taking into account that  $I-\sum\limits_{i=1}^n S_iS_i^*$  is a rank one operator, we deduce that
\begin{equation*}
\begin{split}
\rank \Delta_T&\leq \rank \left[P_{\cN_J}\left( I-\sum\limits_{i=1}^n S_iS_i^*\right)|\cN_J\otimes I_\cD\right]\\
&\leq \dim \cD.
\end{split}
\end{equation*}
Now, using Theorem \ref{dil1}, we conclude that
$
\text{\rm dil-ind}\,(T)=\rank \Delta_T.
$
\end{proof}

\begin{theorem}\label{dil2} Let $J\neq F_n^\infty$
be a WOT-closed two-sided ideal of $F_n^\infty$ generated by a family
 $\cP_J$  of homogenous polynomials, $\cH$ be a separable Hilbert space, and
 $T:=[T_1,\ldots, T_n]$, \ $T_i\in B(\cH)$, be a row contraction such that
\begin{equation}\label{constr}
p(T_1,\ldots, T_n)=0 \quad \text{ for any }\quad p\in \cP_J.
\end{equation}
Then there exists a separable Hilbert space $\cK_\pi$ and a $*$-representation
$\pi:C^*(B_1,\ldots, B_n)\to B(\cK_\pi)$  which annihilates the compact operators and
$$
\pi(B_1)\pi(B_1)^*+\cdots +\pi(B_n)\pi(B_n)^*=I_{ \cK_\pi},
$$
such that
\begin{enumerate}
\item[(i)]
$\cH$ can be identified with a $*$-cyclic co-invariant subspace of
$\tilde\cK:=(\cN_J\otimes \overline{\Delta_T\cH})\oplus \cK_\pi$
under the operators
$$
V_i:=\left[\begin{matrix}
B_i\otimes I_{\overline{\Delta_T\cH}}&0\\0&\pi(B_i)
\end{matrix}\right],\quad i=1,\ldots,n;
$$
\item[(ii)]
$
T_i^*=V_i^*|\cH,\quad i=1,\ldots, n.
$
\end{enumerate}
  \end{theorem}
\begin{proof}
According to \cite{Po-unitary},
if $\cP_J$ consists of homogeneous polynomials, then
$$\text{\rm range}\,K_{T,r}\subseteq \cN_J\otimes \cH\quad \text{
for any } \ r\in(0,1),
$$
the constrained Poisson kernel $K_{J,T,r}:=(P_{\cN_J}\otimes I_\cH)K_{T,r}$ is an isometry,
and there is a unique unital completely contractive linear map
$$\Phi_{J,T}:\overline{\text{\rm span}}\{B_\alpha B_\beta^*:\ \alpha,\beta\in \FF_n^+\}\to B(\cH)
$$
such that $\Phi_{J,T}(B_\alpha B_\beta^*)=T_\alpha T_\beta^*$, \ $\alpha,\beta\in \FF_n^+$.
Applying Arveson extension theorem \cite{Ar1}
to the map $\Phi_{J,T}$, we obtain a unital
completely positive linear map $\Psi_{J,T}:C^*(B_1,\ldots, B_n)\to B(\cH)$.
Let $\tilde\pi:C^*(B_1,\ldots, B_n)\to B(\tilde\cK)$ be a minimal Stinespring dilation
of $\Psi_{J,T}$, i.e.,
$$\Psi_{J,T}(X)=P_{\cH} \tilde\pi(X)|\cH,\quad X\in C^*(B_1,\ldots, B_n),
$$
and $\tilde\cK=\overline{\text{\rm span}}\{\tilde\pi(X)h:\ h\in \cH\}.$
Notice that, for each $i=1,\ldots, n$,
\begin{equation*}
\begin{split}
\Psi_{J,T}(B_i B_i^*)&=T_iT_i^*=P_\cH\tilde\pi(B_i)\tilde\pi(B_i^*)|\cH\\
&=P_\cH \tilde\pi(B_i)(P_\cH+P_{\cH^\perp})\tilde\pi(B_i^*)|\cH\\
&=\Psi_{J,T}(B_i B_i^*)+
(P_\cH \tilde\pi(B_i)|_{\cH^\perp})(P_{\cH^\perp} \tilde\pi(B_i^*)|\cH).
\end{split}
\end{equation*}
Hence, we infer that $P_\cH \tilde\pi(B_i)|_{\cH^\perp}=0$ and
\begin{equation}\label{morph}
\begin{split}
\Psi_{J,T}(B_\alpha X)&=P_\cH(\tilde\pi(B_\alpha) \tilde\pi(X))|\cH\\
&=(P_\cH\tilde\pi(B_\alpha)|\cH)(P_\cH\tilde\pi(X)
|\cH)\\
&=\Psi_{J,T}(B_\alpha)
\Psi_{J,T}(X)\end{split}
\end{equation}
for any $X\in C^*(B_1,\ldots, B_n)$ and $\alpha\in \FF_n^+$.
Note that the Hilbert space $\tilde\cK$ is separable and $\cH$ is an invariant subspace under each $\tilde\pi(B_i)^*$, \ $i=1,\ldots, n$, due to the fact that $P_\cH \tilde\pi(B_i)|_{\cH^\perp}=0$. This means that \begin{equation}\label{coiso}
\tilde\pi(B_i)^*|\cH=\Psi_{J,T}(B_i^*)=T_i^*,\quad i=1,\ldots, n.
\end{equation}

Now, since the subspace $\cN_J$ contains the constants, we can apply Theorem
\ref{compact} and deduce that
all the compact operators $\cL\cC(\cN_J))$ in $B(\cN_J)$ are contained in
$C^*(B_1,\ldots, B_n)$.
According to Theorem \ref{wold}, the representation $\tilde\pi$ decomposes into
a direct sum $\tilde\pi=\pi_0\oplus \pi$ on $\tilde \cK=\cK_0\oplus \cK_\pi$,
where $\pi_0$, $\pi$  are disjoint representations of $C^*(B_1,\ldots, B_n)$
on the Hilbert spaces $\cK_0$ and $\cK_\pi$, respectively, such that
\begin{equation}\label{sime}
\cK_0\simeq\cN_J\otimes \cG, \quad  \pi_0(X)=X\otimes I_\cG, \quad  X\in C^*(B_1,\ldots, B_n),
\end{equation}
 for some Hilbert space $\cG$, and $\pi$ is a representation such that
$\pi(\cL\cC(\cN_J))=0$.
Since $P_\CC^{\cN_J}=I-\sum\limits_{i=1}^n B_iB_i^*$ is a rank one projection in
$C^*(B_1,\ldots, B_n)$, we deduce that
$\sum\limits_{i=1}^n \pi(B_i)\pi(B_i)^*=I_{\cK_\pi}$, and
$$
\dim \cG=\dim (\text{\rm range}\,\tilde\pi(P_\CC^{\cN_J})).
$$
Using the minimality of  the Stinespring representation $\tilde\pi$ and the proof of Theorem \ref{compact}, we have
\begin{equation*}\begin{split}
\text{\rm range}\,\tilde\pi(P_\CC^{\cN_J})&=
\overline{\text{\rm span}}\{\tilde\pi(P_\CC^{\cN_J})\tilde\pi(X)h:\ X\in C^*(B_1,\ldots, B_n), h\in \cH\}\\
&=
\overline{\text{\rm span}}\{\tilde\pi(P_\CC^{\cN_J})\tilde\pi(Y)h:\ Y\in \cC_0, h\in \cH\}\\
&=
\overline{\text{\rm span}}\{\tilde\pi(P_\CC^{\cN_J})\tilde\pi(B_\alpha P_\CC^{\cN_J} B_\beta^*)h:\ \alpha,\beta\in \FF_n^+, h\in \cH\}\\
&=
\overline{\text{\rm span}}\{\tilde\pi(P_\CC^{\cN_J})\tilde\pi(B_\beta^*)h:\ \beta\in \FF_n^+, h\in \cH\}.
\end{split}
\end{equation*}
On the other hand,  using  relation \eqref{morph}, we have
\begin{equation*}\begin{split}
\left<\tilde\pi(P_\CC^{\cN_J})\tilde\pi(B_\alpha^*)h,
\tilde\pi(P_\CC^{\cN_J})\tilde\pi(B_\beta^*)k\right>
&=
\left<h,\pi(B_\alpha)\pi(P_\CC^{\cN_J})\pi(B_\beta^*)h\right>\\
&=
\left<h,T_\alpha\left(I_\cK-\sum_{i=1}^n T_iT_i^*\right)T_\beta^*h\right>\\
&=
\left<\Delta_TT_\alpha^*h,\Delta_TT_\beta^*k\right>
\end{split}
\end{equation*}
for any $h, k \in \cH$.
This shows that one can define  a unitary operator $\Lambda:\text{\rm range}\,\tilde\pi(P_\CC^{\cN_J})\to \overline{\Delta_T\cH}$ by setting
$$
\Lambda(\tilde\pi(P_\CC^{\cN_J})\tilde\pi(B_\alpha^*)h):=\Delta_T T_\alpha^*,\quad h\in \cH,
$$
and extending it by linearity and continuity.
Therefore,  we have
$$
\dim[\text{\rm range}\,\pi(P_\CC^{\cN_J})]=
\dim \overline{\Delta_T\cH}=\dim \cG.
$$
Hence,
making the appropriate identification of $\cG$ with $\overline{\Delta_T\cH}$ and using  relations \eqref{coiso} and \eqref{sime},
  we  obtain the  required dilation.
 This completes the proof.
\end{proof}

\begin{corollary}Let $V:=[V_1,\ldots, V_n]$ be the dilation  of Theorem $\ref{dil2}$. Then,
\begin{enumerate}
\item[(i)]
 $V$ is a constrained shift if and only if $T$ is a pure  constrained row contraction;
\item[(ii)]
 $V$ is a Cuntz type  representation if  and only if
$T$ is a  constrained row contraction such that
$$
T_1T_1^*+\cdots +T_nT_n^*=I.
$$
\end{enumerate}
\end{corollary}
\begin{proof}
Notice that
$$
\sum_{|\alpha|=k}T_\alpha T_\alpha^*=P_\cH
\left[\begin{matrix}
\sum\limits_{|\alpha|=k}B_\alpha B_\alpha^*\otimes I_{\overline{\Delta_T\cH}}&0\\0& I_{\cK_\pi}\end{matrix}\right]|\cH
$$
and therefore,
$$
\text{\rm SOT-}\sum_{|\alpha|=k}T_\alpha T_\alpha^*=P_\cH
\left[\begin{matrix}
 0&0\\0& I_{\cK_\pi}\end{matrix}\right]|\cH.
$$
This shows that
$T$ is a pure row contraction if and only if $P_\cH P_{\cK_\pi} |\cH=0$. The latter condition is equivalent to $\cH\perp (0\oplus \cK_\pi)$, which implies $\cH\subset
\cN_J\otimes\overline{\Delta_T\cH}$.
Now, since $\cN_J\otimes\overline{\Delta_T\cH}$ is reducing for each operator $V_i$, \ $i=1,\ldots, n$,
and $\tilde \cK$ is the smallest reducing subspace  for the same operators, which contains $\cH$, we conclude that
$\tilde\cK=\cN_J\otimes\overline{\Delta_T\cH}$, which proves (i).

Now assume that the dilation  $V$ is a Cuntz type  representation, i.e., $\sum\limits_{i=1}^n V_iV_i^*=I_{\tilde\cK}$. Since
$$\sum_{|\alpha|=k} V_\alpha V_\alpha^*=
\left[\begin{matrix}
\sum\limits_{|\alpha|=k}B_\alpha B_\alpha^*\otimes I_{\overline{\Delta_T\cH}}&0\\0& I_{\cK_\pi}\end{matrix}\right],
$$
we deduce that
$\sum\limits_{|\alpha|=k}B_\alpha B_\alpha^*\otimes I_{\overline{\Delta_T\cH}}=I_{\cK_0}$
for any $k=1,2,\ldots$. Due to the fact that  $\text{\rm SOT-}\lim\limits_{k\to\infty}
\sum\limits_{|\alpha|=k}B_\alpha B_\alpha^*=0$, we must have
$\cK_0=\{0\}$. Using the proof of Theorem \ref{dil2}, we get $\cG=\{0\}$, which means $\Delta_T=0$.
The proof is complete.
\end{proof}

Under additional hypotheses, one can obtain the following  remarkable particular case of Theorem $\ref{dil2}$, where the dilation is unique up to a unitary equivalence.

\begin{corollary}\label{unique}
If in addition to the hypotheses of Theorem $\ref{dil2}$
\begin{equation}\label{BB^*}
\overline{\text{\rm span}}\,\{B_\alpha B_\beta^*:\ \alpha,\beta\in \FF_n^+\}=C^*(B_1,\ldots, B_n),
\end{equation}
then the  dilation of Theorem $\ref{dil2}$ is
minimal, i.e., $\tilde\cK=\bigvee\limits_{\alpha\in \FF_n^+} V_\alpha \cH$,   and  it is unique up to a unitary equivalence.

Let $T':=[T_1',\ldots, T_n']$, \ $T_i'\in B(\cH')$, be  another row contraction subject to the same constraints as $T$ and let $V':=[V_1',\ldots, V_n']$ be the corresponding dilation.  Then $T$ and $T'$ are unitarily equivalent if and only if
$$\dim \overline{\Delta_T\cH}=\dim \overline{\Delta_{T'}\cH'}
$$
and there are unitary operators
 $\Lambda:\cN_J\otimes\overline{\Delta_T\cH}\to \cN_J\otimes\overline{\Delta_{T'}\cH'}$ and
$\Gamma:\cK_\pi\to \cK_{\pi '}$ such that
$$
\Lambda (B_i\otimes I_{\overline{\Delta_T\cH}})=(B_i\otimes I_{\overline{\Delta_{T'}\cH'}})\Lambda\ \text{ and }\  \Gamma\pi(B_i)=\pi'(B_i) \Gamma\quad
\text{ for } \quad i=1,\ldots, n,
$$
and
$$\left[ \begin{matrix} \Lambda&0\\
0&\Gamma\end{matrix} \right] \cH=\cH'.
$$
\end{corollary}

\begin{proof} A closer look at the proof of  Theorem \ref{dil2} reveals that, under condition \eqref{BB^*}, the map $\Psi_{J,T}$ is unique.
Using the uniqueness  of the minimal Stinespring representation (see \cite{St},  \cite{Arv}), one can prove the uniqueness of the minimal dilation of Theorem \ref{dil1}. The last part of this corollary follows using  standard  arguments concerning
representation theory  of $C^*$-algebras \cite{Arv-book} and the uniqueness of minimal completely positive dilations of completely positive maps of $C^*$-algebras.
\end{proof}

In what follows we  present   several examples when the condition \eqref{BB^*} is satisfied.

\begin{example} \label{Exe}
Let $\cP_J\subset F_n^\infty$ by a set of polynomials and let $J$ be the WOT-closed two-sided ideal of $F_n^\infty$ generated by $\cP_J$. The condition \eqref{BB^*} is satisfied   in the following particular cases.

\begin{enumerate}
\item[(i)] If $\cP_J:=0$, then $\cN_J=F^2(H_n)$, $B_i=S_i$, and therefore $S_j^*S_i=\delta_{ij} I$. In this case, Theorem $\ref{dil2}$ and Corollary $\ref{unique}$ imply the  standard noncommutative isometric dilation theorem  for row contraction
\cite{Po-isometric}.
\item[(ii)] If
$\cP_J:=\{S_iS_j-S_jS_i:\ i,j=1,\ldots, n\}$, then $\cN_J=F^2_s$, the symmetric Fock space,  and $B_i$, $i=1,\dots, n$, are the creation operators on the symmetric Fock space. We obtain in this case the dilation theorem for commuting row contractions $($see \cite{Dr}, \cite{Arv}, and \cite{Po-poisson}$)$.
\item[(iii)]
If $B_1,\ldots, B_n$ are essentially  normal.
\item[(iv)] Let $\cP_J$ be a set  of homogenous polynomials in $F_n^\infty$.
According to Lemma $\ref{inv.sub}$, $U\cN_J$ is a subspace
invariant under $S_i^*$,\ $i=1,\ldots, n$. Using the characterization of the invariant subspace  for the left creation operators \cite{Po-charact}, there exists an essentially unique sequence $\{\varphi_p(S_1,\ldots, S_n)\}_{p=1}^N\subset F_n^\infty$, \  $N=1,2,\ldots, \infty$, of isometries with orthogonal ranges such that
$$P_{\cN_J}=I-\sum_{p=1}^N \varphi_p(S_1,\ldots, S_n)\varphi_p(S_1,\ldots, S_n)^*,
$$
where the series is SOT-convergent if $N=\infty$.
If the above sequence is finite ($N<\infty$) and $\varphi_p(S_1,\ldots, S_n)$,\ $p=1,\ldots, N$, are in the noncommutative disc algebra $\cA_n$, then condition
\eqref{BB^*} holds. Indeed, in this case we have
$$B_i^* B_j=P_{\cN_J}S_i^*\left(I-\sum_{p=1}^N \varphi_p(S_1,\ldots, S_n)\varphi_p(S_1,\ldots, S_n)^*\right)S_j|\cN_J, \quad i,j=1,\ldots, n.
$$
Since $S_i^*S_j=\delta_{ij}I$ and $\cN_J$ is invariant under $S_i^*$, \ $i=1,\ldots, n$, we deduce that $B_i^*B_j$ is  in\
$\overline{\text{\rm span}}\,\{B_\alpha B_\beta^*:\ \alpha,\beta\in \FF_n^+\}$.
\item[(v)]
If $\cP_J:=\{S_\alpha:\ |\alpha|= m\}$, then $P_{\cN_J}=I-\sum\limits_{|\alpha|=m}S_\alpha S_\alpha^*$. In this case, we have
$$B_i^* B_j=P_{\cN_J}S_i^*\left(I-\sum_{|\alpha|=m}S_\alpha S_\alpha^*\right)S_j|\cN_J
$$
 and $B_i^*B_j$ is  in
$\text{\rm span}\,\{B_\alpha B_\beta^*:\ \alpha,\beta\in \FF_n^+\}$.
\item[(vi)] If $\cP_J:=\{S_\alpha:\ |\alpha|= m\}\cup \{S_iS_j-S_jS_i:\ i,j=1,\ldots, n\}$.
\item[(vii)] If $\cP_J:=\{S_jS_i-q_{ji} S_iS_j: \ i<j, \ i,j=1,\ldots, n\}$ for some $q_{ji}\in \CC$,  then  $B_i^*B_j=0$ if $i\neq j$ and $B_i^*B_i$ can be written as a linear combination of the identity and $B_jB_j^*$, \ $j=1,\ldots,n$. In this case we obtain the  dilation result from
\cite{BB}.
\end{enumerate}
\end{example}

\smallskip

Let $T:=[T_1,\ldots, T_n]$, \ $T_i\in B(\cH)$, be a row contraction and let
$\cC\subset F_n^\infty$. If $T$ is a c.n.c.  row contraction and $\varphi \in \cC$, then  $\varphi(T_1,\ldots, T_n)$ is defined by the $F_n^\infty$-functional calculus
for row contractions \cite{Po-funct}. When $T$ is an arbitrary row contraction, then  we assume that $\cC$ consists of polynomials.

Denote by $\cM_\cC$ the closed span of all co-invariant spaces $\cM\subseteq \cH$
under $T_1,\ldots, T_n$ such that
$$
\varphi(P_{\cM}T_1 |\cM,\ldots, P_{\cM}T_n|\cM)=0\quad \text{ for any }\ \varphi\in \cC.
$$
We call the row contraction
$$[P_{\cM_\cC}T_1 |\cM_\cC,\ldots, P_{\cM_\cC}T_n|\cM_\cC]
$$
the {\it maximal $\cC$-constrained piece} of $[T_1,\ldots, T_n]$.

\begin{lemma}\label{max-con}
If $V:=[V_1,\ldots, V_n]$, \ $V_i\in B(\cH)$, is a row contraction then
\begin{equation*}
\begin{split}
\cM_\cC&=
\overline{\text{\rm span}}\left\{ V_\alpha \varphi(V_1,\ldots, V_n)\cH:\ \varphi\in \cC, \alpha\in \FF_n^+\right\}^\perp \\
&=\bigcap _{\varphi\in \cC, \alpha\in \FF_n^+}
\varphi(V_1,\ldots, V_n)^*V_\alpha^*.
\end{split}
\end{equation*}
\end{lemma}

\begin{proof}
Denote
 $$
\cE:=\text{\rm span}\left\{ V_\alpha \varphi(V_1,\ldots, V_n)\cH:\ \varphi\in \cC, \alpha\in \FF_n^+\right\}
$$
and note that $\cE^\perp$ is co-invariant under $V_1,\ldots, V_n$.
If $h\in \cE^\perp$ and $k\in \cH$, then
$$
0=\left<\varphi(V_1,\ldots, V_n)k, h\right>=\left< k, \varphi(V_1,\ldots, V_n)^*h\right>.
$$
Hence, we get
$\varphi(P_{\cE^\perp}V_1 |\cE^\perp,\ldots, P_{\cE^\perp}V_n|\cE^\perp)=0.$
Let $\cM$ be a co-invariant subspace under $V_1,\ldots, V_n$ such that
$\varphi(P_{\cM}V_1 |\cM,\ldots, P_{\cM}V_n|\cM)=0$.
For any
$h\in \cM$ and $\alpha\in \FF_n^+$, we have  $V_\alpha^*h\in \cM$, therefore $\varphi(V_1,\ldots, V_n)^*V_\alpha^*h=0$. This implies
$
\left<h, V_\alpha \varphi(V_1,\ldots, V_n)k\right>=0$
for any $k\in \cH$, which shows that $\cM\subseteq \cE^\perp$ and completes the proof.
\end{proof}

Using this lemma and the definition of the subspace $\cN_J$, one can easily prove the following.

\begin{proposition}\label{max-piece}
Let $J\neq F_n^\infty$ be an arbitrary WOT-closed two-sided ideal of $F_n^\infty$,
 $S_1,\ldots, S_n$ be the left creation operators on the full Fock space $F^2(H_n)$, and
$B_i:=P_{\cN_J}S_i|\cN_J$, \ $i=1,\ldots, n$.
Then
  $[B_1,\ldots, B_n]$  is the maximal $J$-constrained piece of $[S_1,\ldots, S_n]$.
\end{proposition}

We consider now the particular  case  when the row contraction  $T:=[T_1,\ldots, T_n]$ is pure, i.e.,
$\text{\rm SOT-}\lim\limits_{k\to\infty} \sum\limits_{|\alpha|=k} T_\alpha T_\alpha^*=0$. In  this case, the results of this section can be extended to a larger class of constrained row contractions.

\begin{theorem}\label{dil3}
Let $J\neq F_n^\infty$
be a WOT-closed two-sided ideal of $F_n^\infty$ and
let $[T_1,\ldots, T_n]$, \ $T_i\in B(\cH)$, be a pure row contraction such that
$$
f(T_1,\ldots, T_n)=0\quad \text{ for any }\quad f\in J.
$$
Then the following statements hold:
\begin{enumerate}
\item[(i)]
The constrained Poisson kernel $K_{J,T}:\cH\to \cN_J\otimes \overline{\Delta_T\cH}$ defined by setting
$$K_{J,T}:=(P_{\cN_J}\otimes I)K_T$$
 is an isometry, $K_{J,T}\cH$ is co-invariant under
$B_i\otimes  I_{\overline{\Delta_T\cH}}$, $i=1,\dots, n$,   and
\begin{equation}\label{Pois}
T_i=K_{J,T}^*(B_i\otimes  I_{\overline{\Delta_T\cH}})K_{J,T},\quad i=1,\ldots, n.
\end{equation}
\item[(ii)] $\text{\rm dil-ind}\,(T)=\text{\rm rank}\,\Delta_T$.
\item[(iii)]
If $1\in \cN_J$, then the dilation provided by \eqref{Pois} is minimal.
\end{enumerate}

If,  in addition,  $1\in \cN_J$ and
\begin{equation}\label{C*}
\overline{\text{\rm span}}\{B_\alpha B_\beta^*:\ \alpha,\beta\in \FF_n^+\}=
C^*(B_1,\ldots, B_n),
\end{equation}
then  we have:
\begin{enumerate}
\item[(iv)] The minimal  dilation provided by \eqref{Pois} is unique up to an isomorphism.
\item[(v)] The minimal  dilation  is the maximal $J$-constrained piece
of the standard noncommutative  isometric dilation of $T$.
\item[(vi)]
A pure constrained row contraction has  $\rank \Delta_T=m$,  $m=1,2,\ldots, \infty$, if and only if it is unitarily equivalent to one obtained by compressing
$[B_1\otimes I_{\CC^m},\ldots, B_n\otimes I_{\CC^m}]$ to a co-invariant subspace $\cM\subset \cN_J\otimes \CC^m$  under  $B_1\otimes I_{\CC^m},\ldots, B_n\otimes I_{\CC^m}$, with the property that $\dim P_0\cM=m$, where $P_0$ is the orthogonal projection of $\cN_J\otimes \CC^m$ onto the subspace $1\otimes \CC^m$.
\end{enumerate}
\end{theorem}
 \begin{proof}
Part (i)  follows from \cite{APo}. To prove (ii), let $\cD$ be a Hilbert space  such that $\cH$ can be identified with a co-invariant subspace of $\cN_J\otimes \cD$ under $B_i\otimes I_\cD$, \ $i=1,\ldots, n$, and  such that
 $T_i=P_\cH(B_i\otimes I_\cD)|\cH$ for  $i=1,\ldots, n$.
Then
\begin{equation*}
\begin{split}
I_\cH-\sum_{i=1}^n T_iT_i^*&=
P_\cH^{\cN_J\otimes \cD}\left[ \left(I_{\cN_J}-\sum_{i=1}^n B_iB_i^*\right)\otimes I_\cD
\right]|\cH\\
&=
P_\cH^{\cN_J\otimes \cD}\left[ P_{\cN_J}(I-\sum_{i=1}^n S_iS_i^*)|\cN_J\otimes I_\cD
\right]|\cH\\
&=
P_\cH^{\cN_J\otimes \cD}\left[ P_{\cN_J}P_\CC|\cN_J\otimes I_\cD
\right]|\cH.
\end{split}
\end{equation*}
Hence, $\rank \Delta_T\leq \dim \cD$. Using (i), we deduce that  the dilation index of $T$ is equal to $\rank \Delta_T$.

Assume now that $1\in \cN_J$. As in the proof of Theorem \ref{cyclic}, we
obtain $P_\CC^{\cN_J} P_{\cN_J}e_\alpha=0$ for any $\alpha\in \FF_n^+$, $|\alpha|\geq 1$.
On the other hand, the definition of the constrained Poisson kernel $K_{J,T}$  implies
$$
P_0K_{J,T}h=\lim_{m\to \infty}\sum_{k=0}^m \sum_{|\alpha|=k} P_\CC^{\cN_J} P_{\cN_J} e_\alpha\otimes \Delta_T T_\alpha^*h,\quad h\in \cH,
$$
where $P_0:=P^{\cN_J}_\CC \otimes I_{\overline{\Delta_T \cH}}$. Therefore,
$P_0K_{J,T}\cH=\overline{\Delta_T \cH}$.  Using Corollary \ref{cyc} in the particular case when $\cM:=K_{J,T}\cH$ and $\cD:=\overline{\Delta_T \cH}$, we deduce that
$K_{J,T}\cH$ is cyclic for  $B_i\otimes I_\cD$,\ $i=1,\ldots, n$, which proves  the minimality of the dilation \eqref{Pois}, i.e.,
\begin{equation}\label{minimal1}
\cN_J\otimes \overline{\Delta_T \cH}=\bigvee_{\alpha\in \FF_n^+} (B_\alpha\otimes I_{\overline{\Delta_T \cH}}) K_{J,T}\cH.
\end{equation}

Now we assume  that $1\in\cN_J$ and   that relation \eqref{C*} holds. Consider another minimal dilation of $T$, i.e.,
\begin{equation}
\label{another}
T_i=V^* (B_i\otimes I_\cD)V,
\end{equation}
where $V:\cH\to \cN_J\otimes \cD$ is an isometry, $V\cH$ is co-invariant under
$B_i\otimes I_\cD$,\ $i=1,\ldots, n$, and
\begin{equation}\label{minimal2}
\cN_J\otimes \cD=\bigvee_{\alpha\in \FF_n^+} (B_\alpha\otimes I_{\cD}) V\cH.
\end{equation}
We know  (see \cite{APo}) that  there exists a  unital completely positive linear map
$$\Phi:\text{\rm span}\{B_\alpha B_\beta^*:\ \alpha,\beta\in \FF_n^+\}\to B(\cH)$$ such that
$\Phi(B_\alpha B_\beta^*)=T_\alpha T_\beta^*$, \ $\alpha, \beta\in \FF_n^+$. Due to \eqref{C*}, $\Phi$ has a unique extension to
$C^*(B_1,\ldots, B_n)$. Consider the $*$-representations
\begin{equation*}
\begin{split}
\pi_1:C^*(B_1,\ldots, B_n)\to B(\cN_J&\otimes \overline{\Delta_T \cH}),\quad \pi_1(X):= X\otimes I_{\overline{\Delta_T \cH}}\\
\pi_2:C^*(B_1,\ldots, B_n)\to B(\cN_J&\otimes \cD),\quad \pi_2(X):= X\otimes I_{\cD}.
\end{split}
\end{equation*}
It is easy to see  that due to relations  \eqref{Pois}, \eqref{another}, \eqref{C*}, and the co-invariance of the subspaces $K_{J,T}\cH$ and $V\cH$, we have
$$
\Phi(X)=K_{J,T}^*\pi_1(X)K_{J,T}=V^*\pi_2(X)V,\quad X\in C^*(B_1,\ldots, B_n).
$$
Now, due to the minimality conditions \eqref{minimal1} and  \eqref{minimal2}, and relation \eqref{C*}, we deduce that $\pi_1$ and $\pi_2$ are  minimal Stinespring dilations of $\Phi$.  Since they are unique, there exists a unitary operator $U:\cN_J\otimes \overline{\Delta_T \cH}\to \cN_J\otimes \cD$ such that
\begin{equation}\label{int-U}
U(B_i\otimes I_{\overline{\Delta_T \cH}})=(B_i\otimes I_\cD)U,\quad
i=1,\ldots,n,
\end{equation}
and $UK_{J,T}=V$.  Hence, we also have
$$
U(B_i^*\otimes I_{\overline{\Delta_T \cH}})=(B_i^*\otimes I_\cD)U,\quad
i=1,\ldots,n.
$$
By Theorem \ref{compact},  $C^*(B_1,\ldots, B_n)$ is irreducible, so we must have $U=I_{\cN_J}\otimes W$, where $W\in B(\overline{\Delta_T \cH},\cD)$ is a unitary operator.
Therefore, $\dim \overline{\Delta_T \cH}=\dim\cD$ and
$UK_{J,T}\cH=V\cH$, which proves that the two dilations are unitarily equivalent.

In the particular case when $J=\{0\}$, part (iv) shows that
$$S:=[S_1\otimes I_{\overline{\Delta_T \cH}},\ldots, S_n\otimes I_{\overline{\Delta_T \cH}}]$$
is a realization of the standard minimal isometric dilation of $[T_1,\ldots, T_n]$.
 Using Lemma
\ref{max-con} and Proposition \ref{max-piece}, one can easily see that the maximal $J$-constrained piece of
$S$ coincides with $[B_1\otimes I_{\overline{\Delta_T \cH}},\ldots, B_n\otimes I_{\overline{\Delta_T \cH}}]$, which proves (v).

Now, we prove (vi). The implication ``$\implies$'' follows from part (i).
Conversely, assume that
$$
T_i=P_\cH(B_i\otimes I_{\CC^m})|\cH,\quad i=1,\ldots,n,
$$
where $\cH\subset \cN_J\otimes \CC^m$ is a co-invariant subspace under
$B_i\otimes I_{\CC^m}$,\  $i=1,\ldots,n$,
with $\dim P_0\cH=m$ (recall $P_0:= P_\CC^{\cN_J}\otimes I_{\CC^m}$).
It is clear that $T:=[T_1,\ldots, T_n]$ is a pure $J$-constrained row contraction.
Consider first the case when $m<\infty$. Since $P_0\cH\subseteq \CC^m$ and
$\dim P_0\cH=m$, we deduce that $P_0\cH=\CC^m$. By Corollary \ref{cyc}, we have
\begin{equation}
\label{perpend}
\cH^\perp\cap \CC^m=\{0\}.
\end{equation}
On the other hand, since $I_{\cN_J}-\sum\limits_{i=1}^n B_iB_i^*=P_\CC^{\cN_J}$, we obtain
\begin{equation*}\begin{split}
\rank \Delta_T&=\rank P_\cH\left[\left(I_{\cN_J}-\sum_{i=1}^n B_iB_i^*\right)\otimes I_{\CC^m}\right]|\cH\\
&=\rank P_\cH P_0|\cH=\dim P_\cH P_0\cH\\
&=\dim P_\cH\CC^m.
\end{split}
\end{equation*}
If $\rank \Delta_T<m$, then there exists  a nonzero vector $h\in \CC^m$ with
$P_\cH h=0$, which contradicts  relation \eqref{perpend}. Therefore, we must have $\rank \Delta_T=m$.

Now, we consider the case $m=\infty$. According to Theorem \ref{cyclic}, setting
$\cE:=P_0\cH$, we have
$$
\bigvee_{\alpha\in \FF_n^+} (B_\alpha\otimes I_{\CC^m})\cH=\cN_J\otimes \cE,
$$
which is reducing for $B_i\otimes I_{\CC^m}$, \ $i=1,\ldots, n$.
Therefore,
$$
T_i=P_\cH(B_i\otimes I_\cE)|\cH,\quad i=1,\ldots, n.
$$
Using the uniqueness of the  minimal dilation of $T$ (see (iv)), we deduce that
$$\dim \overline{\Delta_T \cH}=\dim\cE=\infty.
$$
The proof is complete.
\end{proof}

An important question remains open.  Let $J$ be a WOT-closed two-sided ideal of $F_n^\infty$  generated by a family $\cP_J$ of  homogeneous polynomials such that $1\in \cN_J$ and condition \eqref{BB^*} holds. If $[T_1,\ldots, T_n]$ is a arbitrary $J$-constrained row contraction, is it true that the minimal  dilation  provided by Theorem \ref{dil2} is the maximal $J$-constrained piece
of the standard noncommutative  isometric dilation of $T$ $?$ We should mention   that the answer to this question is positive  in the commutative case, i.e., when
$\cP_J:=\{S_iS_j-S_jS_i:\ i,j=1,\ldots, n\}$ (see \cite{BBD}).

Under  certain natural conditions on the ideal $J$,   we can characterize the pure $J$-constrained row contractions of rank one.

\begin{corollary}\label{rank1} Let $J$ be  a WOT-closed two-sided ideal of $F_n^\infty$  such that  $1\in\cN_J$ and condition \eqref{C*} is satisfied.
If $\cM\subset\cN_J$ is  a co-invariant  subspace under $B_1,\ldots B_n$, then the $n$-tuple
$$
T:=[T_1,\ldots, T_n], \quad  T _i:=P_\cM B_i|\cM,\ i=1,\ldots,n,
$$
is a pure $J$-constrained  row contraction such that
$\rank\Delta_T=1$.

If $\cM'$ is another co-invariant subspace for
$B_1,\ldots, B_n$, which gives rise to a row contraction $T'$, then $T$ and $T'$ are unitarily equivalent if and only if $\cM=\cM'$.

Every pure constrained row contraction with $\rank\Delta_T=1$ is unitarily equivalent to one obtained by compressing
$[B_1,\ldots, B_n]$ to a co-invariant subspace  for $B_1,\ldots, B_n$.
\end{corollary}
\begin{proof}
Since $\cM\subset\cN_J$ is  a co-invariant  subspace under  $B_1,\ldots B_n$,
we have
$$
f(T_1,\ldots, T_n)=P_\cM f(B_1,\ldots, B_n)|\cM=0,\quad f\in J,
$$
and
\begin{equation*}\begin{split}
I_\cM-\sum_{i=1}^nT_iT_i^*&=
P_\cM\left(I_{\cN_J}-\sum_{i=1}^n B_iB_i^*\right)|\cM\\
&=P_\cM P_\CC^{\cN_J}|\cM.
\end{split}
\end{equation*}
Hence, $[T_1,\ldots, T_n]$ is a constrained row contraction with $\text{\rm rank}\, \Delta_T\leq 1$. On the other hand, since
$$
\sum_{|\alpha|=k} T_\alpha T_\alpha^*=P_\cM\left(\sum_{|\alpha|=k} B_\alpha B_\alpha^*\right)|\cM, \quad k=1,2,\ldots,
$$
and $[B_1,\ldots, B_n]$ is a pure row contraction, we deduce that
$[T_1,\ldots, T_n]$ is pure. This also implies that $\Delta_T\neq 0$, so
$\text{\rm rank}\, \Delta_T\geq 1$.
Consequently, we have $\rank \Delta_T=1$.

To prove the second part  of this corollary, notice that,  as in the proof of Theorem \ref{dil3} part (iv),  one can show that $T$ and $T'$ are unitarily equivalent if and only if there exists a unitary operator $\Lambda:\cN_J\to \cN_J$ such that
$$
\Lambda B_i=B_i \Lambda, \ \text{ for } i=1,
\ldots, n, \ \text{ and }  \ \Lambda \cM=\cM'.
$$
This implies that $\Lambda$  commutes with $C^*(B_1,\ldots, B_n)$ which, due to Theorem \ref{compact}, is irreducible. Therefore, $\Lambda$ must be  a scalar multiple of the identity. Hence, we have $\cM=\Lambda \cM=\cM'$.

Finally, the last part of this corollary follows from Theorem \ref{dil2}
and Corollay \ref{unique}.
\end{proof}

\begin{remark}
In the particular case when $n=1$,  the noncommutative analytic Toeplitz algebra $F_n^\infty$ can be identified with the classical Hardy algebra
$H^\infty$. Every  $w^*$-closed  ideal of $H^\infty$  has the form $J=\varphi H^\infty$, where $\varphi\in H^\infty$ is a inner function. In this case
$$\cN_J=H^2\ominus \varphi H^2=\ker \varphi(S)^*.
$$

 If $T\in B(\cH)$ is a c.n.u. contraction such that $\varphi(T)=0$, then   $T$ is pure contraction $($see \cite{SzF-book}$)$.  According to Theorem $\ref{dil1}$,
$\cH$ can be identified with a co-invariant subspace of $\cN_J\otimes\overline{\Delta \cH}$ under $S$, and
$$T=P_\cH (B\otimes I_{\overline{\Delta_T\cH}})|\cH,
$$
where $B:=P_{\ker \varphi(S)^*}S| \ker \varphi(S)^*$
and $S$ is the unilateral shift on the Hardy space $H^2$.
\end{remark}

  \bigskip

      \section{Characteristic functions  for  row contractions and factorizations} \label{Characteristic}

The purpose of this section is to provide new properties for the {\it standard characteristic function} associated with an arbitrary row contraction and
show  that
$$I-\Theta_T \Theta_T^*=K_TK_T^*,
$$
where $K_T$ is the   Poisson kernel of $T$.
Consequently,   we will show that the curvature invariant and Euler characteristic
asssociated with a Hilbert module  over $\CC\FF_n^+$ generated by  an arbitrary
 row contraction $T$ can be expressed only in terms of the
 characteristic function of $T$.

The characteristic  function associated with an arbitrary row contraction
$T:=[T_1,\ldots, T_n]$, \ $T_i\in B(\cH)$, was introduce in \cite{Po-charact} (see \cite{SzF-book} for the classical case $n=1$)
and it was proved to be  a complete unitary invariant for
completely non-coisometric (c.n.c.)
row contractions.
Using the characterization of multi-analytic operators on Fock spaces (see \cite{Po-analytic}, \cite{Po-tensor}), one can easily see that the characteristic  function  of $T$
is  a multi-analytic operator
$$
\Theta_T(R_1,\ldots, R_n):F^2(H_n)\otimes \cD_{T^*}\to F^2(H_n)\otimes \cD_T
$$
with the formal Fourier representation
\begin{equation*}
\begin{split}
  -I_{F^2(H_n)}\otimes T+
\left(I_{F^2(H_n)}\otimes \Delta_T\right)&\left(I_{F^2(H_n)\otimes \cH}-\sum_{i=1}^n R_i\otimes T_i^*\right)^{-1}\\
&\left[R_1\otimes I_\cH,\ldots, R_n\otimes I_\cH
\right] \left(I_{F^2(H_n)}\otimes \Delta_{T^*}\right),
\end{split}
\end{equation*}
where $R_1,\ldots, R_n$ are the right creation operators on the full Fock space $F^2(H_n)$.
 Here,  we need to clarify some notations since some of them are different from those considered in \cite{Po-charact}.
The defect operators  associated with a row contraction $T:=[T_1,\ldots, T_n]$
are
\begin{equation*}
\Delta_T:=\left( I_\cH-\sum_{i=1}^n T_iT_i^*\right)^{1/2}\in B(\cH) \quad \text{ and }\quad \Delta_{T^*}:=(I-T^*T)^{1/2}\in B(\cH^{(n)}),
\end{equation*}
while the defect spaces are $\cD_T:=\overline{\Delta_T\cH}$ and
$\cD_{T^*}:=\overline{\Delta_{T^*}\cH^{(n)}}$, where
$\cH^{(n)}$ denotes the direct sum of $n$ copies of $\cH$.
In what follows we need the following result.

\begin{lemma}\label{*-limit}
 If $\Theta(R_1,\ldots, R_n)\in R_n^\infty\bar\otimes B(\cH,\cK)$, then
\begin{equation*}
\text{\rm SOT-}\lim_{r\to 1}\Theta(rR_1,\ldots, rR_n)^*=\Theta(R_1,\ldots, R_n)^*.
\end{equation*}
\end{lemma}

 \begin{proof}
We know that any multi-analytic operator $\Theta(R_1,\ldots, R_n)$ with formal  Fourier representation $$
\Theta(R_1,\ldots, R_n)\sim \sum_{k=0}^\infty\sum_{|\alpha|=k}R_\alpha \otimes \theta_{(\alpha)},\qquad  \theta_{(\alpha)}\in B(\cH,\cK),
$$
 has the  property that
$$
\Theta(R_1,\ldots, R_n)=\text{\rm SOT-}\lim_{r\to 1}\sum_{k=0}^\infty\sum_{|\alpha|=k}r^{|\alpha|}R_\alpha \otimes \theta_{(\alpha)},
$$
where the series converges in the uniform norm   for each $r\in (0,1)$.
Now, note that
for every $\beta\in \FF_n^+$, $h\in \cH$, and $g\in F^2(H_n)\otimes \cK)$, we have
\begin{equation*}
\begin{split}
\left< \Theta(R_1,\ldots, R_n)^*(e_\beta\otimes h),g\right>
&=
\left<e_\beta\otimes h,\Theta(R_1,\ldots, R_n)g\right>\\
&=\left<e_\beta\otimes h, \left(\sum_{\alpha\in \FF_n^+, |\alpha|\leq |\beta|} R_\alpha \otimes \theta_{(\alpha)}\right)
g\right>\\
&=
\left<\left( \sum_{\alpha\in \FF_n^+, |\alpha|\leq |\beta|} R_\alpha^* \otimes \theta_{(\alpha)}^*\right)(e_\beta\otimes h), g\right>.
\end{split}
\end{equation*}
 Therefore,
$$
\Theta(R_1,\ldots, R_n)^*(e_\beta\otimes h)
=\left(\sum_{\alpha\in \FF_n^+, |\alpha|\leq |\beta|} R_\alpha^* \otimes \theta_{(\alpha)}^*\right)(e_\beta\otimes h).
$$
Similarly, we have
$$
\Theta(rR_1,\ldots, rR_n)^*(e_\beta\otimes h)
=\left(\sum_{\alpha\in \FF_n^+, |\alpha|\leq |\beta|}r^{|\alpha|} R_\alpha^* \otimes \theta_{(\alpha)}^*\right)(e_\beta\otimes h).
$$Using the last two equalities, we obtain
\begin{equation*}
\label{lim2}
\lim_{r\to 1}\Theta(rR_1,\ldots, rR_n)^*(e_\beta\otimes h)
=\Theta(R_1,\ldots, R_n)^*(e_\beta\otimes h)
\end{equation*}
for any $\beta\in \FF_n^+$ and  $h\in \cH$.
On the other hand, according to the noncommutative von Neumann inequality,
$$
\|\Theta(rR_1,\ldots, rR_n)^*\|\leq \|\Theta(R_1,\ldots, R_n)^*\| \quad
\text{ for any } \ r\in (0,1).
$$
Hence, and due to the fact that the closed span of all vectors $e_\alpha\otimes h$ with $\beta\in \FF_n^+$, $h\in \cH$, coincides with $F^2(H_n)\otimes \cH$, we deduce (using standard arguments) that
$$
\text{\rm SOT-}\lim_{r\to 1}\Theta(rR_1,\ldots, rR_n)^*=\Theta(R_1,\ldots, R_n)^*.
$$
The proof is complete.
\end{proof}

The following factorization result will play an important role in our investigation.

\begin{theorem}\label{factor}
Let $T:=[T_1,\ldots, T_n]$, \ $T_i\in B(\cH)$,  be a row contraction.
Then
\begin{equation}\label{fa}
I-\Theta_T\Theta_T^*=K_TK_T^*,
\end{equation}
where
$\Theta_T$ is  the characteristic function  of $T$ and $K_T$  is the corresponding Poisson kernel.
\end{theorem}
\begin{proof}

Denoting $\tilde T:=[I_{F^2(H_n)}\otimes T_1,\ldots, I_{F^2(H_n)}\otimes T_n]$ and
$\hat R:=[R_1\otimes I_\cH,\ldots, R_n\otimes I_\cH]$,  the characteristic function of $T$  has the representation
\begin{equation}\label{cha-red}
\Theta_T(R_1,\ldots, R_n) = \text{\rm SOT-}\lim_{r\to 1}
\left[-\tilde T+ \Delta_{\tilde T}\left(I_{F^2(H_n)\otimes \cH} -r\hat R {\tilde T}^*\right)^{-1} r\hat R \Delta_{\tilde T^*}\right].
\end{equation}

Define the operators
$$
A:=\tilde T^*,\ B:=\Delta_{\tilde T^*}, \ C:=\Delta_{\tilde T}, \ D:=-\tilde T, \ \text{ and } \ Z:=r\hat R, \ 0<r<1,.
$$
and notice that
\begin{equation*}
\left(\begin{matrix}
A&B\\C&D
\end{matrix}
\right)
=
\left(\begin{matrix}
\tilde T^*
&\Delta_{\tilde T^*}\\\Delta_{\tilde T}&
-\tilde T
\end{matrix}
\right)
\end{equation*}
is a unitary operator. Therefore,
\begin{equation}\label{unitar}
AA^*+BB^*=I,\ CC^*+DD^*=I, \text{ and } \ AC^*+BD^*=0.
\end{equation}
Define
$$\Phi(Z):=D+C(I-ZA)^{-1} ZB
$$
and notice that using relation \eqref{unitar}, we have
\begin{equation*}
\begin{split}
I-\Phi(Z)\Phi(Z)^*
&=
I-DD^* -C(I-ZA)^{-1}ZBD^*-DB^*Z^*(I-A^*Z^*)^{-1} C^*\\
&
\qquad -C(I-ZA)^{-1}ZBB^*Z^*(I-A^*Z^*)^{-1} C^*\\
&=
CC^*+C(I-ZA)^{-1}ZAC^*+CA^*Z^*(I-A^*Z^*)^{-1} C^*\\
&
\qquad -C(I-ZA)^{-1}ZZ^*(I-A^*Z^*)^{-1} C^*\\
&\qquad
+C(I-ZA)^{-1}ZAA^*Z^*(I-A^*Z^*)^{-1} C^*\\
&=
C(I-ZA)^{-1}\left[(I-ZA)(I-A^*Z^*)+ZA(I-A^*Z^*)\right.\\
& \qquad
\left.+(I-ZA)A^*Z^*-ZZ^*+ZAA^*Z^*\right]
(I-A^*Z^*)^{-1} C^*\\
&=
C(I-ZA)^{-1}(I-ZZ^*)
(I-A^*Z^*)^{-1} C^*.
 \end{split}
\end{equation*}
Therefore,
\begin{equation}\label{PZ}
I-\Phi(Z)\Phi(Z)^*=
C(I-ZA)^{-1}(I-ZZ^*)
(I-A^*Z^*)^{-1} C^*.
\end{equation}

Therefore,  according to our notations, for any $r\in (0,1)$,  the defect operator
$$I-\Theta_T(rR_1,\ldots, rR_n)\Theta_T(rR_1.\ldots, rR_n)^*$$
 is equal to
the product
\begin{equation*}
\begin{split}
\Delta_{\tilde T}&(I-r\hat R\tilde T^*)^{-1}(I-r^2\hat R \hat R^*)(I-r\tilde T \hat R^*)^{-1} \Delta_{\tilde T}\\
&=\left(I\otimes \Delta_T\right)
\left(I-r\sum_{i=1}^n R_i\otimes T_i^*\right)^{-1}
\left[\left(I-r^2\sum_{i=1}^* R_iR_i^*\right)\otimes I\right]\\
&\phantom{xxxxxxxxxxxxxxxxxxxxxxxxxxxxxxxxx}
\left(I-r\sum_{i=1}^n R_i^*\otimes T_i\right)^{-1}
\left(I\otimes \Delta_T\right)\\
&=
\left(\sum_{k=0}^\infty\sum_{|\gamma|=k}r^{|\gamma|} R_\gamma \otimes \Delta_T T_{\tilde \gamma}^*\right)
\left[\left(I-r^2\sum_{i=1}^n R_iR_i^*\right)\otimes I\right]
\left(\sum_{p=0}^\infty\sum_{|\beta|=p}r^{|\beta|} R^*_\beta \otimes  T_{\tilde \beta}\Delta_T\right)\\
&=
\sum_{k,p=0}^\infty \sum_{|\gamma|=k, |\beta|=p}
r^{|\gamma|+|\beta|}R_\gamma \left(I-r^2\sum_{i=1}^n R_iR_i^*\right) R_\beta^* \otimes
\Delta_T T_{\tilde\gamma}^*T_{\tilde\beta}\Delta_T.
\end{split}
\end{equation*}

Now, for every $\alpha,\omega\in \FF_n^+$, \ $h\in \cD_T$, $k\in\cD_{T}$, we have
\begin{equation*}
\begin{split}
\left<[I\right.&\left.-\Theta_T(rR_1,\ldots, rR_n)\Theta_T(rR_1,\ldots, rR_n)^*] (e_\alpha\otimes h), e_\omega \otimes k\right>\\
&=
\sum_{\gamma\in \FF_n^+, |\gamma|\leq |\omega|}
\sum_{\beta\in \FF_n^+, |\beta|\leq |\alpha|}
\left< r^{|\gamma|+|\beta|}R_\gamma \left(I-r^2\sum_{i=1}^n R_iR_i^*\right) R_\beta^* e_\alpha, e_\omega\right>
\left< \Delta_T T_{\tilde\gamma}^*T_{\tilde\beta}\Delta_T h,k\right>
\end{split}
\end{equation*}
Using Lemma \ref{*-limit}, we have
$$\text{\rm SOT-}\lim_{r\to 1} \Theta_T(rR_1,\ldots, rR_n)^*=\Theta_T(R_1,\ldots, R_n)^*.
$$
Therefore,  the above computations imply that

\begin{equation*}\begin{split}
\left<[I-\Theta_T(R_1,\ldots, R_n)\right.&\left.\Theta_T(R_1,\ldots, R_n)^*] (e_\alpha\otimes h), e_\omega \otimes k\right>\\
&=
\sum_{\gamma\in \FF_n^+, |\gamma|\leq |\omega|}
\sum_{\beta\in \FF_n^+, |\beta|\leq |\alpha|}
\left< R_\gamma P_\CC R_\beta^* e_\alpha, e_\omega\right>
\left< \Delta_T T_{\tilde\gamma}^*T_{\tilde\beta}\Delta_T h,k\right>\\
&=
\sum_{\gamma\in \FF_n^+, |\gamma|\leq |\omega|}\left<R_\gamma(1),e_\omega\right>
\left< \Delta_T T_{\tilde\gamma}^*T_{\alpha}\Delta_T h,k\right>\\
&=\left< \Delta_T T_{\omega}^*T_{\alpha}\Delta_T h,k\right>.
\end{split}
\end{equation*}
Here, we used the fact that $P_\CC R_\beta^*e_\alpha \neq 0$ if and only if $\beta=\tilde \alpha$ (recall that $\tilde\alpha$ is the reverse of $\alpha$), and  that $R_\gamma(1)=\tilde \gamma$.
On the other hand,
using the definition of the Poisson kernel associated with a row contraction, we deduce that
$K_T^*(e_\alpha\otimes h)=T_\alpha \Delta_T$
and
\begin{equation*}
\begin{split}
\left<K_TK_T^*(e_\alpha\otimes h), e_\omega\otimes k\right>
&=
\left< K_TT_\alpha \Delta_Th, e_\omega\otimes k\right>\\
&=
\left<\sum_{\gamma\in\FF_n^+} e_\gamma \otimes \Delta_T T_\gamma^* T_\alpha \Delta h,  e_\omega\otimes k\right>\\
&=
\left<  \Delta_T T_\omega^* T_\alpha \Delta_T h,  e_\omega\otimes k\right>
\end{split}
\end{equation*}
for any $h,k\in \cD_T$ and $\alpha,\omega\in \FF_n^+$.
Summing up the above computations,  we deduce that
$$
I-\Theta_T(R_1,\ldots, R_n)\Theta_T(R_1,\ldots, R_n)^*=K_TK_T^*,
$$
which completes the proof.
\end{proof}

We recall that the spectral radius of an $n$-tuple of operators $X:=[X_1,\ldots, X_n]$ is defined by
$$
r(X):=\lim_{k\to\infty}\left\| \sum_{|\alpha|=k} X_\alpha X_\alpha^*\right\|^{1/2k}.
$$

A closer look at the proof of Theorem \ref{factor}  reveals the following factorization result. We should add that the operator $I-\hat X \tilde T^*$ is invertible because
$r(X)<1$.
\begin{corollary}\label{fact2}
Let $T:=[T_1,\ldots, T_n]$, $T_i\in B(\cH)$, be a row contraction  and let $\Theta_T$ be its characteristic function. If $X:=[X_1,\ldots, X_n]$, $X_i\in B(\cK)$, is a row contraction with spectral radius $r(X)<1$, then
$$
I_{\cK\otimes \cD_T}-\Theta_T(X_1,\ldots, X_n)\Theta_T(X_1,\ldots, X_n)^*=
\Delta_{\tilde T}(I-\hat X\tilde T^*)^{-1}(I-\hat X \hat X^*)(I-\tilde T \hat X^*)^{-1} \Delta_{\tilde T},
$$
where $\hat X:=[X_1\otimes I_\cH,\ldots, X_n\otimes I_\cH]$ and  the  other notations  are from the proof of Theorem $\ref{factor}$.
\end{corollary}

Let $\CC \FF_n^+$ be the complex free semigroup algebra  generated by the free semigroup
$\FF_n^+$ with generators $g_1,\dots, g_n$ and neutral element $g_0$.
Any $n$-tuple $T_1,\dots, T_n $ of bounded operators on a  Hilbert space
$\cH$ gives rise to  a Hilbert  (left) module over  $\CC \FF_n^+$
in the natural way
$$
f\cdot h:= f(T_1,\dots, T_n)h, \quad f\in  \CC \FF_n^+, h\in \cH.
$$
  We say that $\cH$ is a contractive $\CC \FF_n^+$-module  if $T:=[T_1,\dots, T_n] $ is a row contraction, which is equivalent to
$$
\|g_1\cdot h_1+\cdots +g_n\cdot h_n\|^2\leq \|h_1\|^2+\cdots +\|h_n\|^2,
\quad
 h_1,\dots, h_n\in \cH.
$$
We say that $\cH$ is of  finite rank if $\rank(\cH):=\rank \Delta_T <\infty$.
The curvature invariant  and Euler characteristic associated with an arbitrary row contraction  $T$ (or the Hilbert module $\cH$ associated with $T$) were introduced and studied
in \cite{Po-curvature} and \cite{Kr}. We recall that
$$
\text{\rm curv}\,(T)=
\lim_{m\to\infty}
\frac
{\text{\rm trace}\,[I-\Phi_T^m(I)]} {1+n+\cdots+n^{m-1}}
$$
and
$$
\chi(T)=
\lim_{m\to\infty}
\frac
{\rank[I-\Phi_T^m(I)]} {1+n+\cdots+n^{m-1}},
$$
where $\Phi_T$ is the completely positive map associated with $T$, i.e., $\Phi_T(X):=\sum\limits_{i=1}^n T_iXT_i^*$.

Using Theorem \ref{factor} and some results from \cite{Po-curvature},
 we  can  show that the curvature  and the Euler characteristic of a row contraction  $T$ can be expressed
only in terms     of the standard  characteristic function  $\Theta_T$.

\begin{theorem}\label{curv-char}
Let $T:=[T_1,\ldots, T_n]$, \ $T_i\in B(\cH)$, be a row contraction  with $\rank \Delta_T<\infty$, and let  $\text{\rm curv}\,(T)$ and $\chi(T)$ denote its curvature and Euler characteristic, respectively.
Then
$$
\text{\rm curv}\,(T)=\rank \Delta_T-\lim_{m\to\infty} \frac
{\text{\rm trace}\,
[\Theta_T\Theta_T^*(P_m\otimes I)]} {n^m}
$$
and
$$
\chi(T)=
\lim_{m\to\infty}
\frac
{\rank[(I-\Theta_T \Theta_T^*)(P_{\leq m}\otimes I)]} {1+n+\cdots+n^{m-1}},
$$
where $P_m$ $($resp. $P_{\leq m}$$)$ is the orthogonal projection of the full Fock space $F^2(H_n)$ onto the subspace of all homogeneous polynomials of degree $m$
$($resp. polynomials of degree $\leq m$$)$.
\end{theorem}

\begin{proof}
According to Theorem 2.3 and Corollary 2.7 from \cite{Po-curvature}, we have
$$
\text{\rm curv}\,(T)=\lim_{m\to\infty}
\frac
{\text{\rm trace}\,[K_T K_T^*(P_m\otimes I)]} {n^m}.
$$
Using  the factorization result of Theorem \ref{factor}, the first result follows.

Now, according to Theorem 4.1 of \cite{Po-curvature}, we have
\begin{equation}\label{chi}
\chi(T)=
\lim_{m\to\infty}
\frac
{\rank[(K_T^*(P_{\leq m}\otimes I)K_T]} {1+n+\cdots+n^{m-1}}.
\end{equation}
Since $K_T^*(P_{\leq m}\otimes I)$ has finite rank, we have
$$\rank[(K_T^*(P_{\leq m}\otimes I)K_T]=\rank [K_T^*(P_{\leq m}\otimes I)].
$$
On the other hand, since $K_T$ is one-to-one  on the range of
$K_T^*(P_{\leq m}\otimes I)$, we also have
$$
\rank [K_T^*(P_{\leq m}\otimes I)]=
\rank [K_T K_T^*(P_{\leq m}\otimes I)].
$$
Hence, using relation  \eqref{chi} and Theorem \ref{factor}, we complete the proof.
\end{proof}

  \bigskip

      \section{Characteristic functions  and models for constrained row contractions}\label{Characteristic-constr}

In this section,
a {\it constrained characteristic function } is associated with any constrained row contraction.
  For {\it pure constrained} row contractions,
we show that
this characteristic function  is  a  complete unitary invariant and provide a   model in terms of it. We  also show that Arveson's curvature invariant and Euler characteristic
asssociated with a Hilbert module over $C[z_1,\ldots, z_n]$ generated by  a
commuting row contraction $T$ can be expressed only in terms of the
constrained characteristic function of $T$.

Let $J$ be a WOT-closed two-sided ideal of the noncommutative analytic Toeplitz algebra $F_n^\infty$ generated by a family of polynomials $\cP_J$.
We define the {\it constrained characteristic  function} associated with a $J$-constrained  row contraction
$T:=[T_1,\ldots, T_n]$, \ $T_i\in B(\cH)$,
to be   the  multi-analytic operator (with respect to the constrained shifts $B_1,\ldots, B_n$)
$$
\Theta_{J,T}(W_1,\ldots, W_n):\cN_J\otimes \cD_{T^*}\to \cN_J\otimes \cD_T
$$
defined by the formal Fourier representation
$$ -I_{\cN_J}\otimes T+
\left(I_{\cN_J}\otimes \Delta_T\right)\left(I_{{\cN_J}\otimes \cH}-\sum_{i=1}^n W_i\otimes T_i^*\right)^{-1}\\
\left[W_1\otimes I_\cH,\ldots, W_n\otimes I_\cH
\right] \left(I_{\cN_J}\otimes \Delta_{T^*}\right).
$$
Taking into account that $\cN_J$ is  a co-invariant
subspace under $R_1,\ldots, R_n$, we can see that  $\Theta_{J,T}$ is the {\it  maximal $J$-constrained piece} of the standard  characteristic function  $\Theta_T$ of the row contraction $T$. More precisely, we have
\begin{equation}\label{rest}
\begin{split}
\Theta_{T}(R_1,\ldots, R_n)^*(\cN_J\otimes\cD_T)
&\subseteq \cN_J\otimes \cD_{T^*}\ \text{  and  }\\
P_{\cN_J\otimes \cD_{T}}\Theta_{T}(R_1,\ldots, R_n)|\cN_J\otimes \cD_{T^*}&=\Theta_{J,T}(W_1,\ldots, W_n).
\end{split}
\end{equation}

Let us remark that   the above definition
of the constrained characteristic function  makes sense (and has  the same properties) when $J$ is an arbitrary WOT-closed two-sided ideal of $F_n^\infty$ and $T:=[T_1,\ldots, T_n]$ is an arbitrary c.n.c.  $J$-constrained row contraction.

\begin{theorem}\label{J-factor}
Let $J\neq F_n^\infty$ be a WOT-closed two-sided ideal of   $F_n^\infty$ generated by a family of polynomials $\cP_J$.
Let $T:=[T_1,\ldots, T_n]$, \ $T_i\in B(\cH)$,  be a $J$-constrained  row contraction.
Then
\begin{equation}\label{J-fa}
I_{\cN_J\otimes \cD_T}-\Theta_{J,T}\Theta_{J,T}^*=K_{J,T}K_{J,T}^*,
\end{equation}
where $\Theta_{J,T}$ is  the  constrained characteristic function
of   $T$ and $K_{J,T}$ is  the corresponding constrained Poisson
kernel.
\end{theorem}

\begin{proof}
 The constrained Poisson kernel associated with $T$ is  $K_{J,T}:\cH \to \cN_J\otimes \overline{\Delta_T\cH}$ defined by
\begin{equation}\label{J-K}
K_{J,T}:=(P_{\cN_J}\otimes I_{\overline{\Delta_T\cH}})K_T,
\end{equation}
where $K_T$ is the standard Poisson kernel of $T$.
According to  the proof of Theorem \ref{dil1} $\text{\rm range}\, K_T\subseteq \cN_J\otimes
\overline{\Delta_T\cH}$.
Using Theorem \ref{factor} and taking the compression of relation
\eqref{fa} to the subspace $\cN_J\otimes \cD_T\subset F^2(H_n)\otimes \cD_T$, we obtain
$$
I_{\cN_J\otimes \cD_T}-P_{\cN_J\otimes \cD_T}\Theta_T(R_1,\ldots, R_n)
\Theta_T(R_1,\ldots, R_n)^*|\cN_J\otimes \cD_T=
P_{\cN_J\otimes \cD_T}K_TK_T^*|\cN_J\otimes \cD_T.
$$
Taking into account  relations
 \eqref{rest},  \eqref{J-K},  and that $W_i^*=R_i^*|\cN_J$, \ $i=1,\ldots, n$,  we infer that
$$
I_{\cN_J\otimes \cD_T}-\Theta_{J,T}(W_1,\ldots,W_n)\Theta_{J,T}(W_1,\ldots,W_n)^*=K_{J,T}K_{J,T}^*.
$$
\end{proof}
As in the proof of Theorem \ref{J-factor}, one can use Corollary \ref{fact2} to obtain the following constrained version of it.

\begin{corollary}\label{fact3} Let $J\neq F_n^\infty$ be a WOT-closed two-sided ideal of   $F_n^\infty$ generated by a family of polynomials $\cP_J$ and
let $T:=[T_1,\ldots, T_n]$, \ $T_i\in B(\cH)$,  be a $J$-constrained  row contraction.
 If $X:=[X_1,\ldots, X_n]$, $X_i\in B(\cK)$, is a $J$-constrained   row contraction with spectral radius $r(X)<1$, then
$$
I_{\cK\otimes \cD_T}-\Theta_{J,T}(X_1,\ldots, X_n)\Theta_{J,T}(X_1,\ldots, X_n)^*=
\Delta_{\tilde T}(I-\hat X\tilde T^*)^{-1}(I-\hat X \hat X^*)(I-\tilde T \hat X^*)^{-1} \Delta_{\tilde T},
$$
where $\hat X:=[X_1\otimes I_\cH,\ldots, X_n\otimes I_\cH]$ and  the  other notations  are from the proof of Theorem $\ref{factor}$.
\end{corollary}

Now we present a model for  pure constrained row contractions  in terms of characteristic functions.

\begin{theorem}\label{model}
Let $J\neq F_n^\infty$ be a WOT-closed two-sided ideal of $F_n^\infty$ and $T:=[T_1,\ldots, T_n]$ be a  pure $J$-constrained   row contraction.
Then  the constrained characteristic function \linebreak$\Theta_{J,T}\in \cW(W_1,\ldots, W_n)\bar\otimes B(\cD_{T^*},\cD_T)$ is a partial isometry and
$T$ is unitarly equivalent to the   row contraction
\begin{equation}\label{HH}
\left[P_
{\HH_{J,T}} (B_1\otimes I_{\cD_T})|\HH_{J,T},\ldots,
P_{\HH_{J,T}} (B_n\otimes I_{\cD_T})|\HH_{J,T}\right],
\end{equation}
where $P_{\HH_{J,T}}$ is the orthogonal projection of $\cN_J\otimes \cD_T$ on the Hilbert space
$$\HH_{J,T}:=\left(\cN_J\otimes \cD_T\right)\ominus
\Theta_{J,T}(\cN_J\otimes \cD_{T^*}).
$$
 \end{theorem}
\begin{proof}
According to Theorem \ref{dil3}, the constrained Poisson kernel $K_{J,T}:\cH\to \cN_J\otimes \overline{\Delta_T\cH}$
 is an isometry, $K_{J,T}\cH$ is  a co-invariant  subspace under
$B_i\otimes  I_{\overline{\Delta_T\cH}}$, $i=1,\dots, n$,   and
\begin{equation}\label{Pois2}
T_i=K_{J,T}^*(B_i\otimes  I_{\overline{\Delta_T\cH}})K_{J,T},\quad i=1,\ldots, n.
\end{equation}
Consequently, $K_{J,T}K_{J,T}^*$ is the orthogonal projection of $\cN_J\otimes \overline{\Delta_T\cH}$ onto
$K_{J,T}\cH$. According to  Theorem \ref{J-factor}, relation \eqref{J-fa} shows that $K_{J,T}K_{J,T}^*$ and
$\Theta_{J,T}\Theta_{J,T}^*$ are mutually orthogonal projections such that
$$
K_{J,T}K_{J,T}^*+\Theta_{J,T}\Theta_{J,T}^*=
I_{\cN_J\otimes \overline{\Delta_T\cH}}.
$$
Therefore,
$$
K_{J,T}\cH=(\cN_J\otimes \cD_T)\ominus
\Theta_{J,T}(\cN_J\otimes \cD_{T^*}),
$$
Now,  since $K_{J,T}$ is an isometry, we identify $\cH$ with $\HH_{J,T}:=K_{J,T}\cH$ and, using \eqref{Pois2}, we
deduce that $T$ is unitarily equivalent to the row contraction given by \eqref{HH}.
 This completes the proof.
\end{proof}

Let  $\Phi\in \cW(W_1,\ldots, W_n)\bar\otimes B(\cK_1, \cK_2)$
and
$\Phi'\in \cW(W_1,\ldots, W_n)\bar\otimes B(\cK_1', \cK_2')$ be two multi-analytic operators with respect to $B_1,\ldots, B_n$. We say that $\Phi$ and $\Phi'$ coincide
if there are two  unitary multi-analytic operators
$U_j: \cN_J\otimes \cK_j\to \cN_J\otimes\cK_j'$
such that the diagram
\begin{equation*}
\begin{matrix}
\cN_J\otimes \cK_1& \mapright\Phi&
\cN_J\otimes \cK_2\\
 \mapdown{U_1}& & \mapdown{U_2}\\
\cN_J\otimes \cK_1'&  \mapright{\Phi'}&
\cN_J\otimes \cK_2'
\end{matrix}
\end{equation*}
is commutative, i.e., $\Phi' U_1=U_2\Phi$.
Since
$$U_j(B_i\otimes I_{\cK_1})=(B_i\otimes I_{\cK_1'})U_j, \quad i=1,\ldots, n,
$$
and $U_j$ are unitary operators, we also deduce that
$$U_j(B_i^*\otimes I_{\cK_1})=(B_i^*\otimes I_{\cK_1'})U_j, \quad i=1,\ldots, n.
$$
Taking into account that $C^*(B_1,\ldots, B_n)$ is irreducible (see Theorem \ref{compact}), we conclude that
$$
U_j=I_{\cN_J}\otimes \tau_j, \quad j=1,2,
$$
for some unitary operators $\tau_j\in B(\cK_j, \cK_j')$.

The next result shows that the constrained characteristic function is a complete unitary invariant for pure constrained row contractions.

\begin{theorem}\label{u-inv}
Let $J\neq F_n^\infty$ be a WOT-closed two-sided ideal of $F_n^\infty$ and let $T:=[T_1,\ldots, T_n]$, \ $T_i\in B(\cH)$, and  $T':=[T_1',\ldots, T_n']$,\ $T_i'\in B(\cH')$,  be two $J$-constrained  pure row contractions.
Then $T$ and $T'$ are  unitarily equivalent if and only if their  constrained characteristic functions $\Theta_{J,T}$  and
$\Theta_{J,T'}$ coincide.
\end{theorem}
\begin{proof}
Assume that $T$ and $T'$ are unitarily equivalent and let $U:\cH\to \cH'$ be a unitary operator such that
$T_i=U^*T_i'U$ for any $i=1,\ldots, n$. Simple computations reveal that
$$
U\Delta_T=\Delta_{T'}U \quad \text{ and }\quad
(\oplus_{i=1}^n U)\Delta_{T^*}=\Delta_{T'^*}(\oplus_{i=1}^n U).
$$
Define the unitary operators $\tau$ and $\tau'$ by
setting
$$\tau:=U|\cD_T:\cD_T\to \cD_{T'} \quad \text{ and }\quad
\tau':=(\oplus_{i=1}^n U)|\cD_{T^*}:\cD_{T*}\to \cD_{T'^*}.
$$
Taking into account the definition of the constrained characteristic function, it is easy to see that
$$
(I_{\cN_J}\otimes \tau)\Theta_{J,T}=\Theta_{J,T'}(I_{\cN_J}\otimes \tau').
$$

Conversely, assume that the characteristic functions  of $T$ and $T'$ coincide. According to the remarks preceding the theorem, there exist unitary operators
$\tau:\cD_T\to \cD_{T'}$ and $\tau_*:\cD_{T^*}\to \cD_{{T'}^*}$ such that the following diagram
\begin{equation*}
\begin{matrix}
\cN_J\otimes \cD_{T^*}& \mapright{\Phi_{J,T}}&
\cN_J\otimes \cD_T\\
 \mapdown{I_{\cN_J}\otimes \tau_*}& & \mapdown{I_{\cN_J}\otimes \tau}\\
\cN_J\otimes \cD_{{T'}^*}&  \mapright{\Phi_{J,T'}}&
\cN_J\otimes \cD_{T'}
\end{matrix}
\end{equation*}
is commutative, i.e.,
\begin{equation}\label{com}
(I_{\cN_J}\otimes \tau)\Phi_{J,T}=\Phi_{J,T'}(I_{\cN_J}\otimes \tau_*).
\end{equation}
Hence, we deduce that
$$
\Gamma:=(I_{\cN_J}\otimes \tau)|\HH_{J,T}:\HH_{J,T}\to \HH_{J,T'}
$$
is a unitary operator, where $\HH_{J,T}$ and $ \HH_{J,T'}$
are the model spaces for $T$ and $T'$, respectively (see Theorem \ref{model}).
Since
$$
(B_i^*\otimes I_{\cD_T})(I_{\cN_J}\otimes \tau^*)=
(I_{\cN_J}\otimes \tau^*)(B_i^*\otimes I_{\cD_{T'}}),\quad i=1,\ldots, n,
$$
and $\HH_{J,T}$ (resp. $ \HH_{J,T'}$) is a co-invariant subspace under $B_i\otimes I_{\cD_T}$ (resp. $B_i\otimes I_{\cD_{T'}}$), \ $i=1,\ldots, n$, we deduce that
$$
\left[(B_i^*\otimes I_{\cD_T})|\HH_{J,T}\right] \Gamma^*=
\Gamma^*\left[(B_i^*\otimes I_{\cD_{T'}})|\HH_{J,T'}\right],\quad i=1,\ldots, n.
$$
Hence, we obtain
$$
\Gamma\left[P_{\HH_{J,T}}\left(B_i\otimes I_{\cD_T}\right)|\HH_{J,T}\right]=
\left[P_{\HH_{J,T'}}\left(B_i\otimes I_{\cD_{T'}}\right)|\HH_{J,T'}\right]
\Gamma,\quad i=1,\ldots, n.
$$
Now, using Theorem \ref{model}, we conclude that $T$ and $T'$ are unitarily equivalent. The proof is complete.
\end{proof}

\begin{theorem}\label{new-inv}
Let $J$ be a WOT-closed two-sided ideal of $F_n^\infty$ such that $1\in \cN_J$ and condition \eqref{BB^*}  is satisfied. If $\cM\subseteq \cN_J$ is an invariant subspace under $B_1,\ldots, B_n$, and
$$
T:=[T_1,\ldots, T_n],\quad T_i:=P_{\cM^\perp}B_i |\cM^\perp, \quad i=1,\ldots, n,
$$
then
$$
\cM=\Theta_{J,T}\left(\cN_J\otimes \cD_{T^*}\right),
$$
where $\Theta_{J,T}$ is the  constrained characteristic function of $T$.
 \end{theorem}
\begin{proof}
According to Corollary \ref{rank1}, $T$ is a pure $J$-constrained  row contraction with $\rank\Delta_T=1$.
Therefore, we can identify  the subspace $\cD_T$ with $\CC$.
Hence, and due to Theorem \ref{model}, we have
$$
\HH_{J,T}=\cN_J\ominus \Theta_{J,T}\left(\cN_J\otimes \cD_{T^*}\right)
$$ and  $T$ is unitarily equivalent to
$$
\left[P_{\HH_{J,T}}B_1|\HH_{J,T},\ldots,
P_{\HH_{J,T}}B_n|\HH_{J,T}\right].
$$
Using again Corollary \ref{rank1}, we deduce that $\HH_{J,T}=\cM^\perp$ and therefore
$\cM=\Theta_{J,T}\left(\cN_J\otimes \cD_{T^*}\right)$.
This completes the proof.
\end{proof}

\begin{remark}\label{all}  All the results of this section
apply,  in   particular, to  the    cases discussed
in  Example $\ref{Exe}$.
\end{remark}

{\bf The commutative case.}  As mentioned in Example \ref{Exe}, if
$$
\cP_{J_c}:=\{S_iS_j-S_jS_i:\ i,j=1,\ldots, n\},
$$
then $\cN_{J_c}=F^2_s$, the symmetric Fock space,  and $B_i:=P_{F_s^2} S_i|_{F_s^2}$, $i=1,\dots, n$, are the creation operators on the symmetric Fock space.
Arveson showed in  \cite{Arv} that $F_s^2$  can be identified with his space $H^2$,
  of analytic functions in $\BB_n$, namely,
the
reproducing kernel Hilbert space
with reproducing kernel $K_n: \BB_n\times \BB_n\to \CC$ defined by
 $$
 K_n(z,w):= {\frac {1}
{1-\langle z, w\rangle_{\CC^n}}}, \qquad z,w\in \BB_n.
$$
 The algebra
$W_n^\infty:=P_{F_s^2} F_n^\infty|_{F_s^2}$,
    was    proved to be  the $w^*$-closed algebra generated by  the operators
   $B_i$, \ $i=1,\dots, n$, and the identity (see \cite{APo}).
     Moreover, Arveson showed in  \cite{Arv}
   that $W_n^\infty$ can be identified with  the algebra of all  multipliers  of $H^2$.  Under this identification the creation operators $B_1,\ldots, B_n$ become the multiplication operators $M_{z_1},\ldots, M_{z_n}$ by the coordinate functions $z_1,\ldots, z_n$ of $\CC^n$.

Let $T:=[T_1,\ldots, T_n]$, \ $T_i\in B(\cH)$,  be a $J_c$-constrained row contraction, i.e.,
$$
T_iT_j=T_jT_i, \quad i,j=1,\ldots, n.
$$
Under the above-mentioned identifications, the constrained characteristic function of $T$  is the multiplier  (multiplication operator)
$\Theta_{J_c,T}:H^2\otimes \cD_{T^*}\to H^2\otimes \cD_T$
defined by    the operator-valued analytic function on the
open unit ball
$$\BB_n:=\{z:=(z_1,\ldots, z_n)\in \CC_n:\ |z|=(|z_1|^2+\cdots +|z_n|^2)^{1/2}<1\},
$$
given by
$$
\Theta_{J_c,T}(z):=
-T+\Delta_T(I-z_1T_1^*-\cdots -z_nT_n^*)^{-1} [z_1I_\cH,\ldots, z_nI_\cH]\Delta_{T^*},\quad z\in \BB_n.
$$
 For  sake of simplicity,  we  are going to use the same notation for the multiplication operator $M_{\Theta_{J_c,T}}\in B(H^2\otimes \cD_{T^*},  H^2\otimes \cD_T)$ and its symbol
$\Theta_{J_c,T}$, which is a $B(\cD_{T^*}, \cD_T)$-valued
bounded analytic function  in  $\BB_n$.

All the results  of this section can be written in this commutative setting.

Using Theorem \ref{J-factor}  and Corollary \ref{fact3} (in the commutative case) and some results from \cite{Po-curvature},
 we   show that Arveson's  curvature  and  Euler characteristic associated with  a commutative row contraction  $T$
with  $\rank \Delta_T<\infty$ can be expressed in terms     of the constrained  characteristic function  $\Theta_{J_c,T}$.

\begin{theorem}\label{Arv}
Let $T:=[T_1,\ldots, T_n]$, \ $T_i\in B(\cH)$, be a commutative  row contraction  with $\rank \Delta_T<\infty$, and let  $K(T)$ and $\chi(T)$ denote Arveson's curvature and Euler characteristic, respectively.
Then
\begin{equation*}\begin{split}
K(T)&=\int_{\partial \BB_n}\lim_{r\to 1}\text{\rm trace}\, [I-\Theta_{J_c,T}(r\xi)\Theta_{J_c,T}(r\xi)^*] d\sigma(\xi)\\
&=
\rank \Delta_T-(n-1)! \lim_{m\to \infty}
\frac {\text{\rm trace}\, [\Theta_{J_c,T}\Theta_{J_c,T}^* (Q_m\otimes I_{\cD_T})]}
{n^m},
\end{split}\end{equation*}
where $Q_m$ is the projection of $H^2$ onto the subspace   of homogeneous polynomials of degree $m$,
and
$$
\chi(T)=n!
\lim_{m\to \infty}
\frac {\rank\left[\left(I-\Theta_{J_c,T}\Theta_{J_c,T}^*\right) \left(Q_{\leq m}\otimes I_{\cD_T}\right)\right]} {m^n},
$$
where $Q_{\leq m}$ is the projection of $H^2$ onto the subspace  of all polynomials of degree $\leq m$.
\end{theorem}

\begin{proof}
Using the factorization  result of Corollary \ref{fact3}
in our particular case, we obtain
 \begin{equation*}\begin{split}
I-\Theta_{J_c,T}(z)&\Theta_{J_c,T}(z)^*\\
&=
(1-|z|^2)\Delta_T(I-z_1T_1^*-\cdots -z_nT_n^*)^{-1}
(I-\overline{z}_1T_1-\cdots -\overline{z}_nT_n)^{-1}
\Delta_T
\end{split}
\end{equation*}
for any $z\in \BB_n$.

The first formula follows from the definition of the curvature \cite{Arv2} and the above-mentioned
factorization  for the constrained characteristic function of $T$.
Using Corollary 2.4. and Corollary 2.8 from \cite{Po-curvature},  we have
$$
K(T)=(n-1)!\lim_{m\to \infty}
\frac {\text{\rm trace}\,[(P_m\otimes I)K_TK_T^* (P_m\otimes I)]} {m^{n-1}},
$$
where  $K_T$ is the Poisson kernel of $T$ and $P_m$ is the orthogonal projection of $F^2(H_n)$ onto the subspace of all homogeneous polynomials of degree $m$. Since $T$ is a commutative row contraction, i.e., $J_c$-constrained,  we have
$\text{\rm range}\,K_T\subset F^2_s\otimes \cD_T$ and the constrained Poisson kernel satisfies the equation
$K_{J_c,T}=(P_{F^2_s}\otimes I)K_T$,
where $F^2_s$ is the symmetric Fock space.
Using the standard properties for the {\it trace} and the above relation, we deduce that
\begin{equation}\label{KJ1}
K(T)=(n-1)!\lim_{m\to \infty}
\frac {\text{\rm trace}\,[K_{J_c,T}K_{J_c,T}^* (Q_m\otimes I)]} {m^{n-1}},
\end{equation}
where $Q_m:=P_{F^2_s} P_m| F^2_s$ is the projection
 of $F^2_s$ onto the subspace of homogeneous polynomials of degree $m$. According to Theorem \ref{J-factor}, we have
\begin{equation}\label{J_c}
I- \Theta_{J_c,T}\Theta_{J_c,T}^* =K_{J_c,T}K_{J_c,T}^*.
\end{equation}
Taking onto account  relations \eqref{KJ1} and \eqref{J_c}, we deduce the second formula for the curvarure. Here, of course, we used Arveson's identification of the symmetric Fock space  $F_s^2$ with
his space $H^2$.

Arveson \cite{Arv2} showed that his Euler characteristic  satisfies the equation
$$
\chi(T)=n!\lim_{m\to\infty}\frac {\rank[I-\Phi_T^{m+1}(I)]} {m^n},
$$
where $\Phi_T$ is the completely positive map associated with $T$. We proved in \cite{Po-curvature} (see Corollary 4.3) that
\begin{equation}\label{EU-Po}
\chi(T)=n!\lim_{m\to\infty}\frac{\rank [K_T^*(P_{\leq m}\otimes I)K_T]} {m^n},
\end{equation}
where $P_{\leq m}$ is the orthogonal projection of $F^2(H_n)$  on the subspace of all polynomials of degree $\leq m$.
Using again that $\text{\rm range}\,K_T\subset F^2_s\otimes \cD_T$ and the constrained Poisson kernel satisfies the equation
$K_{J_c,T}=(P_{F^2_s}\otimes I)K_T$,
 we deduce that
\begin{equation*}
\begin{split}
\rank [K_T^*(P_{\leq m}\otimes I)K_T]
&=
\rank [K_T^*(P_{F^2_s}\otimes I)(P_{\leq m}\otimes I)
P_{F^2_s}\otimes I)K_T]\\
&=\rank \left[ K_{J_c,T}^*(Q_{\leq m}\otimes I) K_{J_c,T}\right]\\
&=\rank \left[ K_{J_c,T}^*(Q_{\leq m}\otimes I)\right]\\
&=\rank \left[ K_{J_c,Y}K_{J_c,T}^*(Q_{\leq m}\otimes I)\right],
\end{split}
\end{equation*}
where  $Q_{\leq m}$ is the projection of $F_s^2$ onto the subspace of all polynomials of degree $\leq m$. The last two equalities hold since the operator $K_{J_c,T}^*(Q_{\leq m}\otimes I)$ has finite rank and $K_{J_c,T}$ is one-to-one on the range of $K_{J_c,T}^*(Q_{\leq m}\otimes I)$.
Now, using  relation \eqref{EU-Po}, the above equalities,  and  the factorization \eqref{J_c},  we  obtain the last formula of the theorem.
The proof is complete.
\end{proof}

      \bigskip

      \section{Commutant lifting theorem for constrained row contractions}\label{Commutant}

In this section  we provide a  Sarason \cite{S} type commutant lifting theorem for pure constrained row contractions and   obtain a
Nevanlinna-Pick  \cite{N} interpolation
 result  in our setting.
  Let  $[T_1,\ldots, T_n]$, $T_i\in B(\cH)$,  be a pure row contraction, and let $J$ be a WOT-closed two-sided ideal of $F_n^\infty$  such  that
\begin{equation}\label{J}
\varphi(T_1,\ldots, T_n)=0 \quad \text{ for any  } \ \varphi(S_1,\ldots, S_n)\in J,
\end{equation}
where $\varphi(T_1,\ldots, T_n)$ is  defined using the $F_n^\infty$-functional calculus for  row contractions.
 According to Section \ref{Dilations}, any pure constrained row contraction is unitarily equivalent to the compression of $[B_1\otimes I_\cK,\ldots, B_n\otimes I_\cK]$ to a co-invariant subspace $\cE$ under each operator $B_i\otimes I_\cK$, $i=1,\ldots, n$. Therefore, we have
$$T_i=P_\cE(B_i\otimes I_\cK)|\cE,\quad i=1,\ldots,n.
$$

The following result is  a commutant lifting theorem for pure
constrained row contractions.

\begin{theorem}\label{CLT}
Let $J\neq F_n^\infty$ be a WOT-closed two-sided ideal of  the noncommutative analytic Toeplitz algebra $F_n^\infty$  and let $[B_1,\ldots, B_n]$  and $[W_1,\ldots, W_n]$ be the corresponding constrained shifts acting
on $\cN_J$.  For each $j=1,2$, let $\cK_j$ be a Hilbert space and $\cE_j\subseteq \cN_J\otimes \cK_j$ be a co-invariant subspace  under each operator $B_i\otimes I_{\cK_j}$, \ $i=1,\ldots, n$.
If $X:\cE_1\to \cE_2$ is a bounded operator such that
\begin{equation}\label{int}
X[P_{\cE_1}(B_i\otimes I_{\cK_1})|_{\cE_1}]=[P_{\cE_2}(B_i\otimes I_{\cK_2})]|_{\cE_2}X,\quad i=1,\ldots,n,
\end{equation}
then there exists
$$G(W_1,\ldots, W_n)\in \cW(W_1,\ldots, W_n)\bar\otimes B(\cK_1,\cK_2)$$
such that $G(W_1,\ldots, W_n)^*\cE_2\subseteq \cE_1$,
$$
G(W_1,\ldots, W_n)^*|\cE_2=X^*, \quad \text{ and }\quad \|G(W_1,\ldots, W_n)\|=\|X\|.
$$
In particular, if  $\cE_j:=\cG\otimes \cK_j$, where $\cG $ is  a co-invariant subspace under each  operator $B_i$  and $W_i$, \ $i=1,\ldots, n$, then the above implication becomes an equivalence.
\end{theorem}

\begin{proof}
According to Lemma \ref{inv.sub}, the subspace  $\cN_J\otimes \cK_j$ is  invariant under  each operator $S_i^*\otimes I_{\cK_j}$, \ $i=1,\ldots, n$,  and
$$
(S_i^*\otimes I_{\cK_j})|\cN_J\otimes \cK_j=B_i^*\otimes I_{\cK_j},\quad  i=1,\ldots, n.
$$
Since $\cE_j\subseteq \cN_J\otimes \cK_j$ is invariant   under   $B_i^*\otimes I_{\cK_j}$ it is also invariant under $S_i^*\otimes I_{\cK_j}$  and
$$
(S_i^*\otimes I_{\cK_j})|\cE_j=  (B_i^*\otimes I_{\cK_j})|\cE_j,
\quad i=1,\ldots, n.
$$
Hence, relation \eqref{int} implies
\begin{equation}\label{int2}
XP_{\cE_1}(S_i\otimes I_{\cK_1})|_{\cE_1}=P_{\cE_2}(S_i\otimes I_{\cK_2})|_{\cE_2}X,\quad i=1,\ldots,n.
\end{equation}
For each $j=1,2$, the $n$-tuple $[S_1 \otimes I_{\cK_j},\ldots, S_n \otimes I_{\cK_j}]$ is  an isometric dilation of the row contraction
$$
[P_{\cE_j}(S_1\otimes I_{\cK_j})|\cE_j,\ldots, P_{\cE_j}(S_n\otimes I_{\cK_j})|\cE_j].
$$
 Applying the noncommutative commutant lifting theorem
(\cite{Po-isometric}, \cite{Po-intert}), we find a multi-analytic operator  $\Phi(R_1,\ldots, R_n)\in R_n^\infty\bar\otimes B(\cK_1,\cK_2)$ such that $\Phi(R_1,\ldots, R_n)^*\cE_2\subseteq \cE_1$,
\begin{equation}\label{clt}
\Phi(R_1,\ldots, R_n)^*|\cE_2=X^*,\quad \text{ and } \quad \|\Phi(R_1,\ldots, R_n)\|=\|X\|.
\end{equation}

Let $G(W_1,\ldots, W_n):=P_{\cN_J\otimes \cK_2} \Phi(R_1,\ldots, R_n)|\cN_J\otimes \cK_1$.
According to the remarks preceding Theorem \ref{Beur}, we have
$$
G(W_1,\ldots, W_n)\in [P_{\cN_J}R_n^\infty |\cN_J]\bar \otimes B(\cK_1, \cK_2)
=\cW(W_1,\ldots, W_n)\bar\otimes B(\cK_1, \cK_2).
$$
Since $\Phi(R_1,\ldots, R_n)^*(\cN_J\otimes \cK_2)\subseteq \cN_J\otimes \cK_1$
and $\cE_j\subseteq \cN_J\otimes \cK_j$, relation\eqref{clt} implies
$$
G(W_1,\ldots, W_n)^*\cE_2\subseteq \cE_1\quad \text{ and }\quad G(W_1,\ldots, W_n)^*|\cE_2=X^*.
$$
Hence, and using again \eqref{clt}, we have
$$
\|X\|\leq \|G(W_1,\ldots, W_n)\|\leq \|\Phi(R_1,\ldots, R_n)\|=\|X\|.
$$
Therefore, $\|G(W_1,\ldots, W_n)\|=\|X\|$.

Now, let us prove the last part of the theorem.
The implication ``$\implies$'' is clear from the first part  of the theorem.
For the converse, let $X=P_{\cG \otimes \cK_2} \Psi(W_1,\ldots, W_n)|\cG \otimes \cK_1$, where $\Psi(W_1,\ldots, W_n)\in \cW(W_1,\ldots, W_n)\bar\otimes B(\cK_1, \cK_2)$. Since $B_iW_j=W_jB_i$ for $i,j=1,\ldots, n$,  we have
$$
(B_i^*\otimes I_{\cK_1})\Psi(W_1,\ldots, W_n)^*=\Psi(W_1,\ldots, W_n)^*
(B_i^*\otimes I_{\cK_2}), \quad i=1,\ldots, n.
$$
Now, taking into account that $\cG$ is  an invariant subspace under each of the operators $B_i^*$  and $W_i^*$, \ $i=1,\ldots, n$, we deduce \eqref{int}.
The proof is complete.
\end{proof}

\begin{corollary}
Let $J\neq F_n^\infty$ be a WOT-closed two-sided ideal of  the noncommutative analytic Toeplitz algebra $F_n^\infty$  and let $B_1,\ldots, B_n$  and $W_1,\ldots, W_n$ be the corresponding constrained shifts acting
on $\cN_J$. If $\cK$ is a Hilbert space and $\cG\subseteq \cN_J$ is an invariant subspace under each operator   $B_i^*$  and $W_i^*$, \ $i=1,\ldots, n$, then
$$
\left\{\left[P_\cG \cW(B_1,\ldots, B_n)|\cG\right]\otimes I_\cK\right\}^\prime=
 [P_\cG \cW(W_1,\ldots, W_n)|\cG]\bar\otimes B(\cK).
$$

\end{corollary}

We remark that Theorem \ref{CLT} can be extended to the following  more general setting. The proof follows exactly the same lines so we shall omit it.

For each $j=1,2$, let
$J_j$ be a  WOT-closed two-sided ideal of $F_n^\infty$  and let $[B_1^{(j)},\ldots, B_n^{(j)}]$    be the corresponding constrained shift acting
on $\cN_{J_j}$.
Let $\cE_j\subseteq \cN_{J_j}\otimes \cK_j$ be an invariant subspace  under each operator ${B_i^{(j)}}^*\otimes I_{\cK_j}$, \ $i=1,\ldots, n$, where $\cK_j$ is a Hilbert space.
If $X:\cE_1\to \cE_2$ is a bounded operator such that
$$
XP_{\cE_1}(B_i^{(1)}\otimes I_{\cK_1})|_{\cE_1}=P_{\cE_2}(B_i^{(2)}\otimes I_{\cK_2})|_{\cE_2}X,\quad i=1,\ldots,n,
$$
then there exists $G\in [P_{\cN_{J_2}} R_n^\infty|\cN_{J_1}]\bar\otimes B(\cK_1,\cK_2)$ such that
$$
P_{\cE_2}G|\cE_1=X \quad \text{ and }\quad \|G\|=\|X\|.
$$

Now we can obtain  the following  Nevanlinna-Pick  interpolation
 result  in our setting.  We only sketch the proof which is similar to  that of Theorem 2.4 from \cite{APo} but uses Theorem \ref{CLT}, and point out  what is  new.

\begin{theorem}\label{Nev}
Let $J$ be  a WOT-closed two-sided ideal of $F_n^\infty$  and let $B_1,\ldots, B_n$    be the corresponding constrained shifts acting
on $\cN_J$. Let $\lambda_1,\ldots, \lambda_k$ be $k$ distinct points in the zero set
$$
\cZ_J:=\{\lambda\in \BB_n:\ f(\lambda)=0 \text{ for any } f\in J\},
$$
and let $A_1,\ldots, A_k\in B(\cK)$.
Then there exists $\Phi(B_1,\ldots, B_n)\in \cW(B_1,\ldots, B_m)\bar\otimes B(\cK)$
such that
$$\|\Phi(B_1,\ldots, B_n)\|\leq 1\quad \text{and }\quad
\Phi(\lambda_j)=A_j,\ j=1,\ldots, k,
$$
if and only if the operator matrix
\begin{equation}\label{semipo}
\left[\frac {I_\cK-A_i A_j^*} {1-\left< \lambda_i,\lambda_j\right>}\right]_{k\times k}
\end{equation}
is positive semidefinite.
\end{theorem}

\begin{proof}

  Let $\lambda_j:=(\lambda_{j1},\dots,\lambda_{jn})\in \CC^n$, $j=1,\dots,k$, and denote $\lambda_{j\alpha}:=\lambda_{ji_1}\lambda_{ji_2}\dots
\lambda_{ji_m}$ if
$\alpha=g_{i_1}g_{i_2}\dots g_{i_m}\in \FF_n^+$,
and $\lambda_{jg_0}:=1$.
Define
$$
z_{\lambda_j}:=\sum_{\alpha\in \FF_n^+} \overline {\lambda}_{j\alpha} e_\alpha,\quad
j=1,2,\dots,k.
$$
Notice that, for any $f\in J$, $\lambda\in \cZ_J$, and $\alpha,\beta\in \FF_n^+$,  we have
$$\left< [S_\alpha f(S_1,\ldots, S_n) S_\beta](1), z_\lambda \right>= {\lambda}_\alpha
 {f(\lambda)}  {\lambda}_\beta=0,
$$
which implies  $z_\lambda\in \cN_J$ for any $\lambda\in \cZ_J$. Note also that, since $B_i^*=S_i^*|\cN_J$ for $i=1,\ldots, n$, we have
$$
B_i^*z_{\lambda_j}=\overline{\lambda}_{ji} z_{\lambda_j}\quad
\text{for }\ i=1,\ldots, n, \text{ and } j=1,\ldots, k.
$$
Define the subspace
$$
\cM:=\text{span} \{ z_{\lambda_j}:\ j=1,\dots,k\}
$$
  and   the operators $X_i\in B(\cM\otimes\cK)$ by setting  $X_i=P_\cM B_i|_\cM\otimes I_\cK$, $i=1,\ldots, n$.
    Since $z_{\lambda_1},\dots, z_{\lambda_k}$ are linearly independent, we can define an operator
    $T\in  B(\cM\otimes\cK)$ by setting
    $$
    T^*(z_{\lambda_j}\otimes h)=
z_{\lambda_j}\otimes A_j^* h
    $$
    for any $h\in \cK$ and   $j=1,\dots, k$.
    Notice that  $TX_i=X_iT$ for  $i=1,\dots, n $.

      Since $\cM$ is  a co-invariant subspace  under  each operator $B_i$, $ i=1,\dots,n$,  we can apply  Theorem
\ref{CLT}
       and find
$\Phi(W_1,\dots, W_n)\in \cW(W_1,\ldots, W_n)\bar\otimes B(\cK)$ such that
\begin{equation}\label{PHI*}
\Phi(W_1,\dots, W_n)^* \cM\subset \cM,\quad
\Phi(W_1,\dots, W_n)^*|\cM=T^*,
\end{equation}
 and  $ \|\Phi(W_1,\dots, W_n)\|=\|T\|$.
As in \cite{APo}, one can prove that
$\Phi(\lambda_j)=A_j$,\ $j=1,\ldots, k$, if and only if
\eqref{PHI*} holds.
Moreover, $ \|\Phi(W_1,\dots, W_n)\|\leq 1$ if and only if
$TT^*\leq I_\cM$, which is equivalent to the fact that
the operator matrix \eqref{semipo}
is  positive semidefinite. This completes the proof.
\end{proof}

We should remark that  in the commutative case
when  $J= J_c$ (see Example \ref{Exe} part (ii)),   we recover the result obtained in \cite{APo}, \cite{Po-interpo}, and \cite{DP}.

       %

      \end{document}